\documentclass[final]{amsart}

\usepackage{color}
\usepackage{graphicx,epstopdf}
\usepackage{amsmath}
\usepackage{amssymb}
\usepackage{amsfonts}
\usepackage{amsthm}
\usepackage[color]{showkeys}

\theoremstyle{remark}
\newtheorem{remark}{Remark}[section]
 
\begin{document}
 \title[]{Numerical study of the transverse stability of the 
 Peregrine solution}

\author{Christian Klein}
\address{Institut de Math\'ematiques de Bourgogne, UMR 5584\\
                Universit\'e de Bourgogne-Franche-Comt\'e, 9 avenue Alain Savary, 21078 Dijon
                Cedex, France\\
    E-mail Christian.Klein@u-bourgogne.fr}

\author{Nikola Stoilov}
\address{Institut de Math\'ematiques de Bourgogne, UMR 5584\\
                Universit\'e de Bourgogne-Franche-Comt\'e, 9 avenue Alain Savary, 21078 Dijon
                Cedex, France\\
    E-mail Nikola.Stoilov@u-bourgogne.fr}
\date{\today}

\begin{abstract}
We generalise 
a previously published approach based on a multi-domain spectral method 
on the whole real line in two ways: firstly, a fully 
explicit 4th order method for the time integration, based on 
a splitting scheme and an implicit Runge--Kutta method for the linear 
part, is presented. Secondly, the 1D code is combined with a Fourier spectral method 
in the transverse variable both for elliptic and hyperbolic NLS 
equations. As an example we study the transverse stability of the 
Peregrine solution, an exact solution to the one dimensional 
nonlinear Schr\"odinger (NLS) equation and thus a $y$-independent 
solution to the 2D NLS. 
 It is shown that the 
Peregrine solution is unstable against all standard perturbations, 
and that some perturbations can even lead to a blow-up for the 
elliptic NLS equation.  
\end{abstract}

 
\thanks{This work was partially supported by 
the ANR-FWF project ANuI - ANR-17-CE40-0035, the isite BFC project 
NAANoD, the ANR-17-EURE-0002 EIPHI and by the 
European Union Horizon 2020 research and innovation program under the 
Marie Sklodowska-Curie RISE 2017 grant agreement no. 778010 IPaDEGAN. 
We thank D. Clamond and C. Kharif for interesting us in the question 
of perturbations of the Peregrine solution in two dimensions}
\maketitle

\section{Introduction}
Nonlinear Schr\"odinger (NLS) equations, see \cite{sulem} for general 
references,  of the form
\begin{equation}
    i\partial_{t}u+u_{xx}+\kappa u_{yy}+V(|u|^{2},x)u=0
    \label{NLS},
\end{equation}
where $u=u(x,t):\mathbb{R}\times\mathbb{R}^{d}\mapsto \mathbb{C}$, 
$\kappa=\pm1$, 
are omnipresent  as asymptotic models for 
modulation of waves in applications such as hydrodynamics, nonlinear optics and 
Bose-Einstein condensates. Here we concentrate on the cubic NLS 
($V=2|u|^{2}$) in $d=2$ dimensions in the elliptic  focusing case 
($\kappa=1$) and the hyperbolic case ($\kappa=-1$). Despite the importance of the NLS 
equations in applications, they are notoriously difficult to treat
numerically. This is due to the fact that solutions can have 
zones of rapid modulated oscillations, called dispersive shock waves 
(DSWs), 
and that the focusing elliptic version (\ref{NLS}) is $L^{2}$ 
\emph{critical}. This means that the equation can have an 
$L^{\infty}$ blow-up in finite time for initial data on 
$H^{1}(\mathbb{R}^{2})$ with an $L^{2}$ 
norm larger than the one of the ground state. The reader is referred to 
\cite{MR,sulem} for an asymptotic description
of the blow-up. For numerical approaches for DSWs in this context see 
for instance \cite{etna} and references therein, for blow-up in NLS 
see \cite{sulem} and more recently \cite{KS}, for a numerical study of 
the hyperbolic NLS \cite{KS15}.

These problems are only aggravated if one wants to study solutions such as the 
Peregrine solution  \cite{Peregrine} to the cubic 1D NLS, 
\begin{equation}
    u_{Per} = \left(1-\frac{4(1+4it)}{1+4x^{2}+16t^{2}}\right)e^{2it}
    \label{peregrine},
\end{equation}
see  \cite{dubard} for generalizations thereof, which do not vanish 
asymptotically for $|x|\to\infty$, but which tend algebraically to a constant in 
modulus. Such solutions are often discussed as a model for rogue 
waves: It is claimed that the Peregrine solution has been 
observed experimentally in rogue wave experiments in hydrodynamics 
\cite{cha1,cha2,cha3}, 
in plasma physics \cite{bailung} and in nonlinear optics 
\cite{kibler}. However, recently it has been shown numerically in 
\cite{KBAA,BK,KH} that the Peregrine solution is unstable, and 
its orbital instability in the 1D NLS has been rigorously proven by Mun\~oz 
\cite{munoz}. Interesting remarks on the Peregrine solution can be 
also found in \cite{Li}. For alternative approaches via NLS equations 
to rogue wave phenomena see for instance \cite{Abl} and 
\cite{BPSS,DBU}. 

For  rapidly decreasing or periodic initial data
the standard numerical approach is to use Fourier methods, see \cite{etna} and 
references. Initial data with compact support allow the application of 
perfectly matched layers \cite{Berenger,Zhengpml} or transparent boundary 
conditions \cite{Zhengtbc} (based on \cite{BFS}) to avoid unwanted reflections at the 
boundaries of the computational domains. The latter techniques, with some modification, can be successfully used in the case of algebraically decreasing data as in (\ref{peregrine}). For example, in \cite{BK} one of the authors presents a multi-domain spectral approach with a 
compactification of the infinite interval, i.e., where infinity is a regular
point on the grid. The use of spectral methods for Schr\"odinger 
equations is recommended since the latter 
are purely dispersive equations. This 
means that the introduction of numerical dissipation should be strongly 
limited if dispersive effects like rapid oscillations are of interest. 

In order to study the transverse stability of the Peregrine solution, 
in this paper we extend the approach of \cite{BK} to the 2D 
NLS for perturbations which are periodic or rapidly decreasing in 
$y$. To this end the 1D approach in \cite{BK} is changed in two ways: 
firstly, we use a Fourier spectral method in $y$ together with the original
approach in 1D from \cite{BK} which gives a spectral method both in 
$x$ and $y$, adapted to solutions with an algebraic decay in $x$ to a 
finite value at infinity, and rapidly decreasing or periodic in 
$y$. Secondly, we replace the implicit fourth order Runge-Kutta (IRK4) method 
for the time integration in \cite{BK} with a fourth order splitting 
where the IRK4 method  will be only applied to the 
linear part of the equation, which can be explicitly solved since it 
is a linear relation. This leads to a fully explicit method the accuracy 
of which is studied at the example of the Peregrine solution.

The Peregrine solution can be seen as an exact $y$-independent, and 
thus not localized in $y$, solution to the 2D NLS 
equation (\ref{NLS}). With the code presented 
in this paper, we are able to study transverse perturbations of the 
Peregrine solution in 2D. We consider various types of such 
perturbations, both localised and non-localised in $x$, and this for 
the elliptic and the hyperbolic NLS. It is shown 
that the Peregrine solution is unstable against all studied 
perturbations, and that in certain cases even a blow-up can be observed.

The paper is organized as follows: in Section 2 we present the 
numerical approach for the time integration for the 1D code and test 
its convergence for the example of the Peregrine solution; the code 
is then extended to 2D. In 
Section 3 we consider perturbations of the Peregrine solution 
localized in $x$ for the elliptic NLS equation. Perturbations not localized in $x$ are studied in 
Section 4. Similar perturbations are considered for the hyperbolic 
NLS in section 5. We add some concluding remarks in section 6. 

\section{Numerical approach}
In this section we summarise the numerical approaches applied in this 
paper. Firstly, we recall the main features of the 1D code \cite{BK}, 
then we present a fully explicit variant for the time 
integration scheme of the code and test it at the example of the 
Peregrine solution. The resulting code is generalised to the 2D NLS 
(\ref{NLS}) where the dependence on the 
transverse variable $y$ is treated with Fourier methods.

\subsection{Spatial discretisation}

Spectral methods, see for instance \cite{trefethen} and references 
therein, are an attractive choice for the numerical solution 
of partial differential equations (PDEs), in particular in higher 
dimensions, because of their excellent approximation properties for smooth 
functions. If the function to be approximated is analytic on the 
considered interval, the numerical error whilst applying spectral 
methods is known to decrease exponentially with the resolution, a 
feature which is called spectral convergence. 
An interesting additional advantage of spectral methods is that they 
minimize the introduction of numerical dissipation in the 
approximated PDEs. This is especially attractive in the present 
context of dispersive PDEs where dissipation could suppress the 
dispersive effects one wants to study as DSWs.

The basic idea of spectral methods is to approximate a function on a 
given interval by a set of global functions on this interval. The 
best known such method works for a periodic setting  and is the
truncated Fourier series, i.e., a trigonometric polynomial. It is well 
known that the numerical error in approximating an analytic periodic 
function via a truncated Fourier series decreases exponentially with the 
number of terms in this sum. In a finite precision approach the same 
is true for rapidly decreasing functions if the period is chosen 
large enough that the function and its first derivatives vanish 
with numerical precision at the domain boundaries. This is the 
approach we will use for the variable $y$.

Nevertheless, Fourier techniques are also known to be much less efficient in 
the case of discontinuous periodic functions. In this case the 
numerical error in approximating the function via a truncated Fourier 
series is known to decrease only linearly with the number of terms in the 
sum, and there is a Gibbs phenomenon at the discontinuity in this 
case. Thus, it is considerably better to approximate functions which 
are only smooth on intervals by a sum of polynomials there. Slowly decreasing functions on the whole real 
line can be spectrally treated through a compactification of the external domain, namely
\begin{equation}
    x\in x_{0}[-1,1],\quad s = 1/x\in[-1,1]\frac{1}{x_{0}},
    \label{2int}
\end{equation}
where $x_{0}>0$ is some constant. 
This means we cover the real line with two domains which are both 
mapped to the interval $[-1,1]$. 

\begin{remark}
    Note that in such a 
\emph{multi-domain} method, in principle an arbitrary number of 
intervals can be used. This allows to allocate resolution were 
needed, whereas a compactification of the real line to a simple 
domain, for instance as in \cite{GO} or via the well known map $s=\arctan(x)$, imposes 
the distribution of the numerical resolution for a given 
discretisation. If infinity is not a regular point as here, 
one can also use two infinite intervals as in \cite{BK}, 
$]-\infty,-x_{0}]$ and $[x_{0},\infty[$. In the context of finite element 
methods, such intervals are called \emph{infinite elements}. 
Compactification techniques are standard for elliptic equations, for 
instance in astrophysics, see \cite{lorene} for a similar approach as 
here.
\end{remark}

For the case of the 1D NLS, this means that we solve in interval I 
($x\in x_{0}[-1,1]$) the equation
\begin{equation}
    iu^{I}_{t}+u^{I}_{xx}+2|u^{I}|^{2}u^{I}=0
    \label{NLSI},
\end{equation}
and in domain II ($s\in[-1,1]/x_{0}$) 
\begin{equation}
    iu^{II}_{t}+s^{4}u^{II}_{ss}+2s^{3}u^{II}_{s}+2|u^{II}|^{2}u^{II}=0
    \label{NLSII}.
\end{equation}
Put together domain I and II cover the whole real line, and the 
solution $(u^{I},u^{II})$ has to be at least  $C^{1}(\mathbb{R})$ and thus 
in particular to be differentiable at the domain boundaries. This 
implies the conditions
\begin{equation}
    u^{I}(\pm x_{0}) = u^{II}(\pm1/x_{0}),\quad u^{I}_{x}(\pm x_{0})= 
    -u^{II}_{s}(\pm1/x_{0})/x_{0}^{2}.
    \label{match}
\end{equation}

On each of the two domains we  introduce standard Chebyshev 
collocation points, $x_{n}=x_{0}\cos(n\pi/N_{I})$, 
$n=0,1,\ldots,N_{I}$, and $s_{n}
=\cos(n\pi/N_{II})/x_{0}$, $n=0,1,\ldots,N_{II}$,  
$N_{I},N_{II}\in\mathbb{N}$ and approximate the considered function by the 
Lagrange interpolation polynomial sampled at the collocation points. The 
derivatives of the function $u$ are approximated by the derivatives 
of the Lagrange polynomial which leads to the action of a 
\emph{Chebyshev differentiation matrix} $D$ on the function, see the 
discussion in \cite{trefethen,WR}. This means that the PDEs 
(\ref{NLSI}) and (\ref{NLSII}) are approximated by the finite 
dimensional system of ordinary differential equations (ODEs)
\begin{align}
    u^{I}_{t}&= iD_{x}^{2}u^{I}+2i|u^{I}|^{2}u^{I},
    \nonumber\\
    u^{II}_{t}& =iD_{s}^{2}u^{II}+2i|u^{II}|^{2}u^{II}
    \label{ODEsys},
\end{align}
where in an abuse of notation we have used the same symbols for the 
functions $u^{I}(x)$ and the vector with the components 
$u^{I}(x_{n})$, $n=0,1,\ldots,N$, and similarly for $u^{II}$. Note 
that we always use an even number of points  $N_{II}+1$ in domain II in 
order to make sure that infinity is not a grid point. This simplifies 
the treatment of the singularity in the equation for $u^{II}$. 

This approach can be generalised to two spatial dimensions in the 
following way, where we assume that $u(x,y)$ is either periodic or 
rapidly decreasing in $y$: in this case, after the discretisation in 
$x$, the components of
vector $u^{I}_{n}$ depend on $y$ (the case for $u^{II}$ is completely 
analogous). The $y$-dependence is then approximated via a truncated 
Fourier series for $y\in L_{y}[-\pi,\pi]$,
$$u^{I}_{n}(y)\approx \sum_{k=-M/2+1}^{M/2}\hat{u}^{I}_{nm}\exp(2\pi 
i k 
y/L_{y}/M).$$
This means we discretise also $y$, $y_{1}, \ldots, y_{M}$ and use a 
\emph{Fast Fourier transform} (FFT) to compute the discrete Fourier 
transform of $u^{I}_{n}$. Derivatives in $y$ are then approximated in 
standard way by multiplication of the Fourier coefficients in the sum by 
a factor $ik$. This yields a spectral method in both $x$ and $y$, 
i.e., a combined Chebyshev-Fourier approach for the spatial 
coordinates. Since for a function both analytic in $x$ and $y$, the 
Chebyshev and the Fourier coefficients have to decrease 
exponentially, we can control the spatial resolution of the solution 
via the behavior of the \emph{spectral coefficients}. Note that the Chebyshev 
coefficients can be computed with a \emph{Fast Cosine Transform} 
which is close to the FFT, see \cite{trefethen}.

\subsection{Explicit time integration scheme}
The spatial discretisation of the previous subsection leaves us with a 
finite dimensional ODE system of the form 
\begin{equation}
    U_{t} = \mathcal{L}U+\mathcal{N}[U]
    \label{sysgen},
\end{equation}
where $U$ is an $(N_{I}+N_{II}+2)\times M$ dimensional matrix, $\mathcal{L}$ is a 
linear operator formed essentially by a discretised version of the 
Laplace or d'Alembert operator in (\ref{NLS}), and where $\mathcal{N}[U]$ is a 
nonlinear operator free of derivatives. Since there are second 
derivatives in the linear part of NLS equations, the resulting system 
is \emph{stiff} which means that explicit time integration schemes 
will be inefficient since stability conditions would impose severe 
restrictions on the size of the time steps. This is 
especially true since the operator $D_{x}^{2}$ has a conditioning of 
order $N_{I}^{4}$, see for instance the discussion in 
\cite{trefethen}, and similarly for $s^{4}D_{s}^{2}+2s^{2}D_{s}$.

We therefore chose in \cite{BK} an implicit time integration scheme  to address these stability issues. 
Concretely, we apply a fourth  order Runge-Kutta (IRK4) scheme, the Hammer-Hollingsworth method, a 2-stage 
Gauss scheme. The general formulation of an $s$-stage Runge--Kutta method 
for the initial value problem
$y'=f(y,t),\,\,\,\,y(t_0)=y_0$ is as follows:
\begin{eqnarray}
 y_{n+1} = y_{n} + h      \underset{i=1}{\overset{s}{\sum}} \, 
 b_{i}K_{i}, \\
 K_{i} = f\left(t_{n}+c_ {i}h,\,y_{n}+h  
 \underset{j=1}{\overset{s}{\sum}} \, a_{ij}K_{j}\right),
 \label{K}
\end{eqnarray}
where $b_i,\,a_{ij},\,\,i,j=1,...,s$ are real numbers and
$c_i=   \underset{j=1}{\overset{s}{\sum}} \, a_{ij}$.  
For the IRK4 method used here, one has
$c_{1}=\frac{1}{2}-\frac{\sqrt{3}}{6}$, 
$c_{2}=\frac{1}{2}+\frac{\sqrt{3}}{6}$, $a_{11}=a_{22}=1/4$,
$a_{12}=\frac{1}{4}-\frac{\sqrt{3}}{6}$, 
$a_{21}=\frac{1}{4}+\frac{\sqrt{3}}{6}$ and $b_{1}=b_{2}=1/2$. 

The system 
following from (\ref{K}) for (\ref{sysgen}) is written in the form 
\begin{align}
    (\hat{1}-ha_{11}\mathcal{L})K_{1} 
    & = \mathcal{L}U(t_{n})
    +ha_{12}\mathcal{L}K_{2}
    +\mathcal{N}\left[U(t_{n})+h\sum_{j=1}^{2}a_{1j}K_{j}\right]
    \nonumber,\\
    (\hat{1}-ha_{22}\mathcal{L})K_{2} & 
    =\mathcal{L}U(t_{n})
    +ha_{21}\mathcal{L}K_{1}
    +\mathcal{N}\left[U(t_{n})+h\sum_{j=1}^{2}a_{2j}K_{j}\right]
    \label{Ksys}.
\end{align}
The matching conditions (\ref{match}) are 
implemented via a $\tau$-method, see \cite{BK}. The resulting 
implicit system for $K_{1}$ and $K_{2}$ requires an iterative 
solution. For this we used a 
simplified Newton method to invert the operators of the left hand side of 
(\ref{Ksys}) and get the new $K_{1}$ and $K_{2}$). The 
iteration converges rapidly in simple cases, but can take more than 
100 iterations to reach high precision in some of the examples
discussed in the following. 

In the linear case $U_{t}=\mathcal{L}U$, the solution to (\ref{Ksys}) 
can be written in the form 
\begin{equation}
    \mathcal{L}_{+}\mathcal{L}_{-}(K_{1}+K_{2})=\left(\hat{1}-\frac{h}{2}\mathcal{L}+\frac{h^{2}}{12}\mathcal{L}^{2}\right) (K_{1}+K_{2})=
2\mathcal{L}U(t_{n})
    \label{lineq},
\end{equation}

where 
$$\mathcal{L}_{\pm}=\hat{1}-\frac{1}{4}\left(1\pm 
\frac{i}{\sqrt{3}}\right)h\mathcal{L}.$$
An inversion of the operator containing $\mathcal{L}^{2}\propto 
D_{x}^{4}$ is not 
recommended since the conditioning of the matrix $D_{x}^{4}$ would be of the 
order $N_{I}^{8}$. Furthermore, in order to get a unique solution,  one would have to impose a $C^{3}$ 
condition at the domain boundaries which is in itself problematic. 
Therefore we used in \cite{BK} even in the linear case the iterative 
approach based on (\ref{Ksys}). In this paper however, we are solving 
(\ref{lineq}) by inverting first $\mathcal{L}_{+}$ and then 
$\mathcal{L}_{-}$ in (\ref{lineq}), 
each of them 
with a $C^{1}$ condition. Thus we only have to deal with the 
conditioning of the operator $D_{x}^{2}$, but  in contrast to the 
iterative approach we have to do that twice. In practice this implies, as will be 
discussed in more detail below for the Peregrine example, that the 
maximally achievable accuracy is of the order of $10^{-8}$ instead of 
$10^{-11}$ in the iterative solution. Nevertheless this is more than sufficient 
for our purposes, 
and it is significantly more efficient as no iterations are needed. 
Thus we have a fully explicit scheme though we apply an implicit 
scheme in the linear part since the resulting equation can be explicitly solved 
in this case.

Using this approach for the nonlinear system (\ref{Ksys}) can be done through 
splitting, i.e., to split the equation into the linear part
(\ref{lineq}) and the nonlinear part $U_{t}=\mathcal{N}[U]$ which can 
be integrated explicitly. The linear part will then be integrated 
with IRK4 in the form (\ref{lineq}). 
The  motivation for splitting methods comes from the Trotter-Kato formula
\cite{TK,Ka}
\begin{equation}\label{e11}
 {\lim}_{n\rightarrow\infty}\left(e^{-tA/n}e^{-tB/n}\right)^{n}=e^{-t\left(A+B\right)}
\end{equation}
where $A$ and $B$ are certain unbounded linear operators, for details 
see \cite{Ka}.

The idea of these methods for an equation of the form 
$U_{t}=\left(A+B\right)U$ is to write the solution as
\[
U(t)=\exp(c_{1}tA)\exp(d_{1}tB)\exp(c_{2}tA)\exp(d_{2}tB)\cdots\exp(c_{k}tA)\exp(d_{k}tB)u(0)
\]
where $(c_{1},\,\ldots,\, c_{k})$ and $(d_{1},\,\ldots,\, d_{k})$
are sets of real numbers that represent fractional time steps.  
Yoshida \cite{Y} gave an approach 
which produces split step methods of any even order. We apply here a 
fourth order splitting. Combined with IRK4 for the linear step, this 
gives a fourth order splitting scheme for NLS with an inexact 
solution of the linear step via an implicit fourth order method. However, 
the latter can be done explicitly because of the linearity of the 
equation. The nonlinear step can be integrated as usual for NLS in 
explicit form since for the equation $iu_{t}+|u|^{2}u=0$, $|u|$ is 
constant in time. Thus we get in total a fully explicit scheme which,  
especially in two dimensions, is much more efficient than the IRK4 scheme 
for the full equation with an iterative solution. 

The numerical accuracy of the time integration can be controlled via 
the conserved quantities of the NLS equation of which there are 3 in 
two dimensions (the 1D cubic NLS is of course completely integrable 
and thus has an infinite number of conserved quantities). Since we 
want to treat solutions as the Peregrine solution which do not tend 
to zero at infinity, but satisfy 
$\lim_{|x|\to\infty}|u|=\lambda$, we consider a combination of the 
energy and the $L^{2}$ norm which is defined in such cases and which 
is a conserved quantity of NLS,
\begin{equation}
    E[u]=\int_{\mathbb{R}^{2}}^{}\left(|u_{x}|^{2}+\kappa|u_{y}|^{2}-|u|^{2}(
    |u|^{2}-\lambda^{2})\right)
    \label{energy}.
\end{equation}
Due to unavoidable numerical errors, this quantity will not be 
exactly conserved in actual computations. However the relative (with 
respect to the initial value) conserved 
$E[u]$ gives an indication of the numerical error (as discussed in 
\cite{etna}, it typically overestimates the numerical accuracy by one 
to two orders of magnitude). Note that (\ref{energy}) is not defined 
for the Peregrine solution on $\mathbb{R}^{2}$. But it is for
$\mathbb{R}\times \mathbb{T}$, the situation we study in this paper. 

\subsection{Test}
As a test of the above code we use the Peregrine solution 
(\ref{peregrine}) which can be seen as a $y$ independent solution of 
the 2D NLS equation. Since this solution does not test the $y$ 
dependence in the code, the results are the same as with the 1D code, 
and  we concentrate on the latter to study the dependence of the 
numerical error on the time step. Concretely we impose the initial 
data $u(x,0)=u_{Per}(x,0)$ and solve the NLS equation for these 
initial data. The numerical and the exact solution are compared for 
$t=1$, as the numerical error $\Delta$ we take the $L^{\infty}$ 
norm of the difference between both solutions. We always use $N_{I}=80$ and 
$N_{II}=75$ collocation points in the respective domains, and vary 
the number $N_{t}$ of time steps between $100$ and $1500$. 

\begin{remark}
    Note that the results in 2D could be different from the presented 
    ones below in 1D due to numerical instabilities, in particular 
    for the focusing elliptic NLS known for a modulational 
    instability. We have tested in detail that this is not the case 
    as long as sufficient resolution in $y$ is provided which is 
    always the case in the examples in the following sections. Thus 
    we could restrict the convergence study in $t$ to the 1D case in 
    order to allocate less computational resources.
\end{remark}

\begin{figure}[htb!]
  \includegraphics[width=0.49\textwidth]{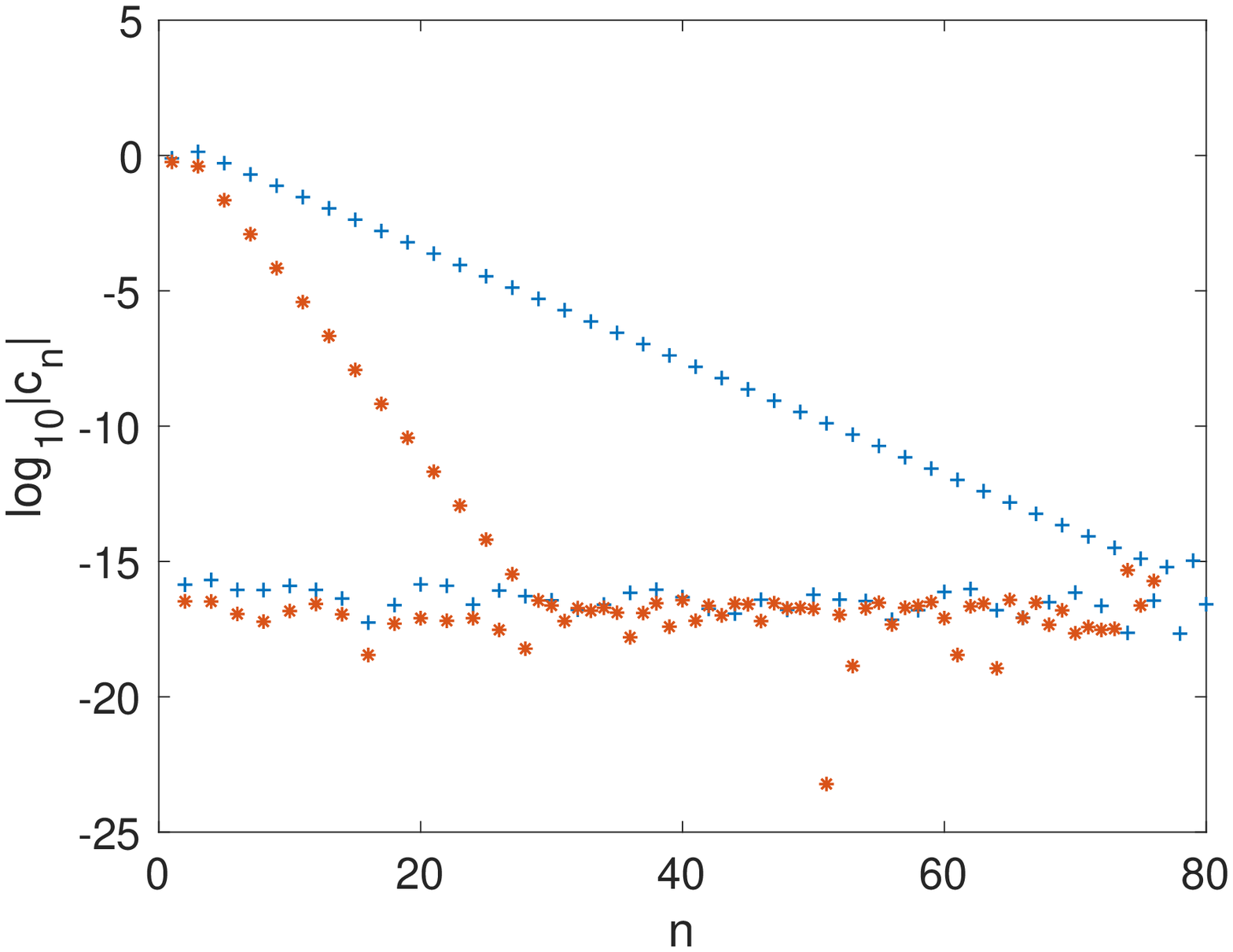}
  \includegraphics[width=0.49\textwidth]{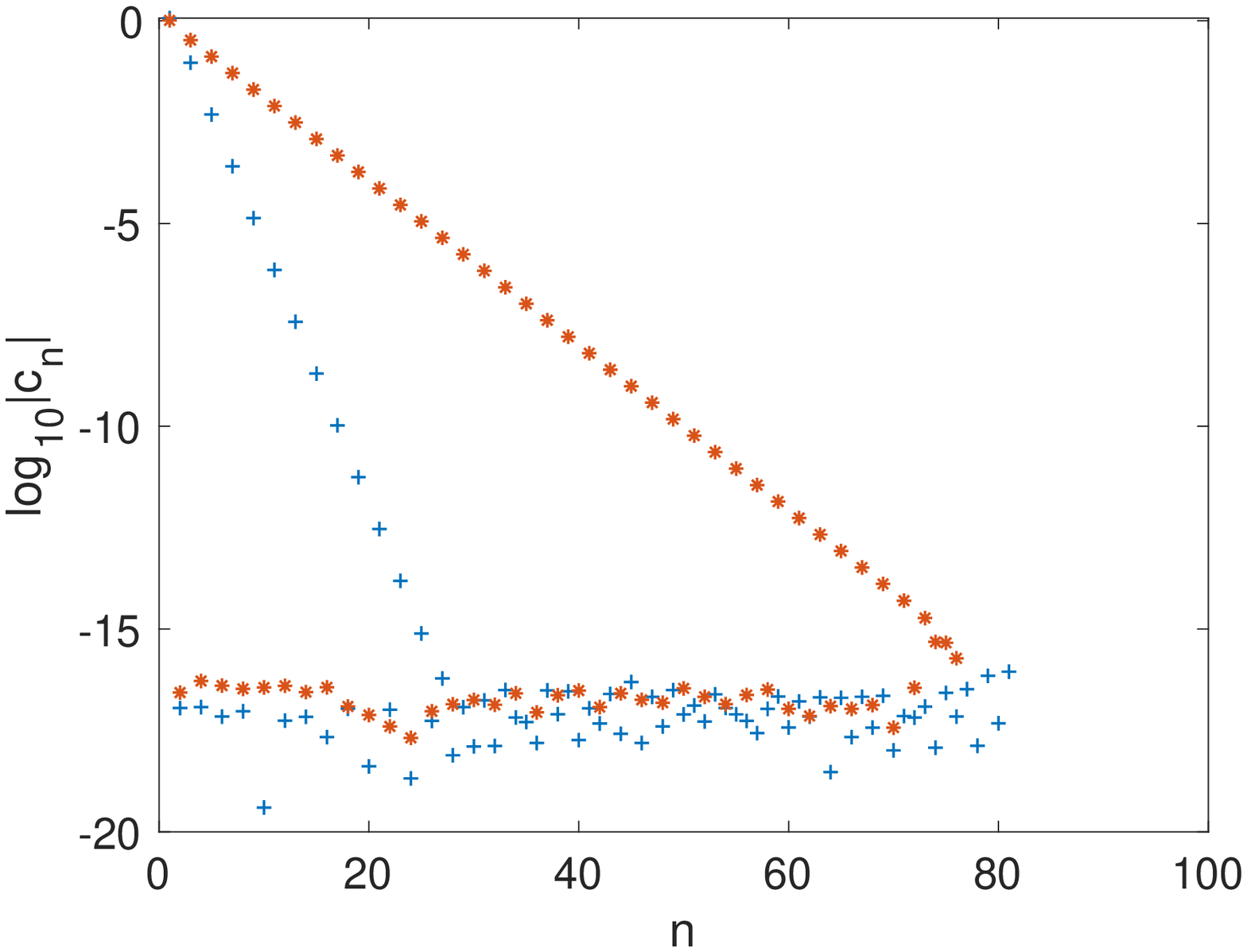}
 \caption{Chebyshev coefficients for the Peregrine solution given by 
 `+' in domain I and `*' in domain II, on the left for $t=0$, on the 
 right for $t=1$.}
 \label{Percoeff}
\end{figure}

The 
Chebyshev coefficients when representing the Peregrine solution in 
the form $u_{Per}(x,0)=\sum_{n=0}^{N}c_{n}T_{n}(x)$ where $T_{n}(x)$, 
$n=0,1,\ldots$ are the Chebyshev polynomials can be seen on the left 
of Fig.~\ref{Percoeff}, the coefficients at the final time on the 
right of the same figure. This shows that the solution is always well 
resolved in space since the Chebyshev coefficients decrease in both 
domains to the order of the rounding error which is of the order of 
$10^{-15}$. Note that for larger values of the number of Chebyshev 
polynomials, the coefficients would just saturate more as in the case 
of the domain I coefficients on the right of figure \ref{Percoeff}, 
and for the coefficients in domain II on the left of the same figure. 

\begin{figure}[htb!]
  \includegraphics[width=0.49\textwidth]{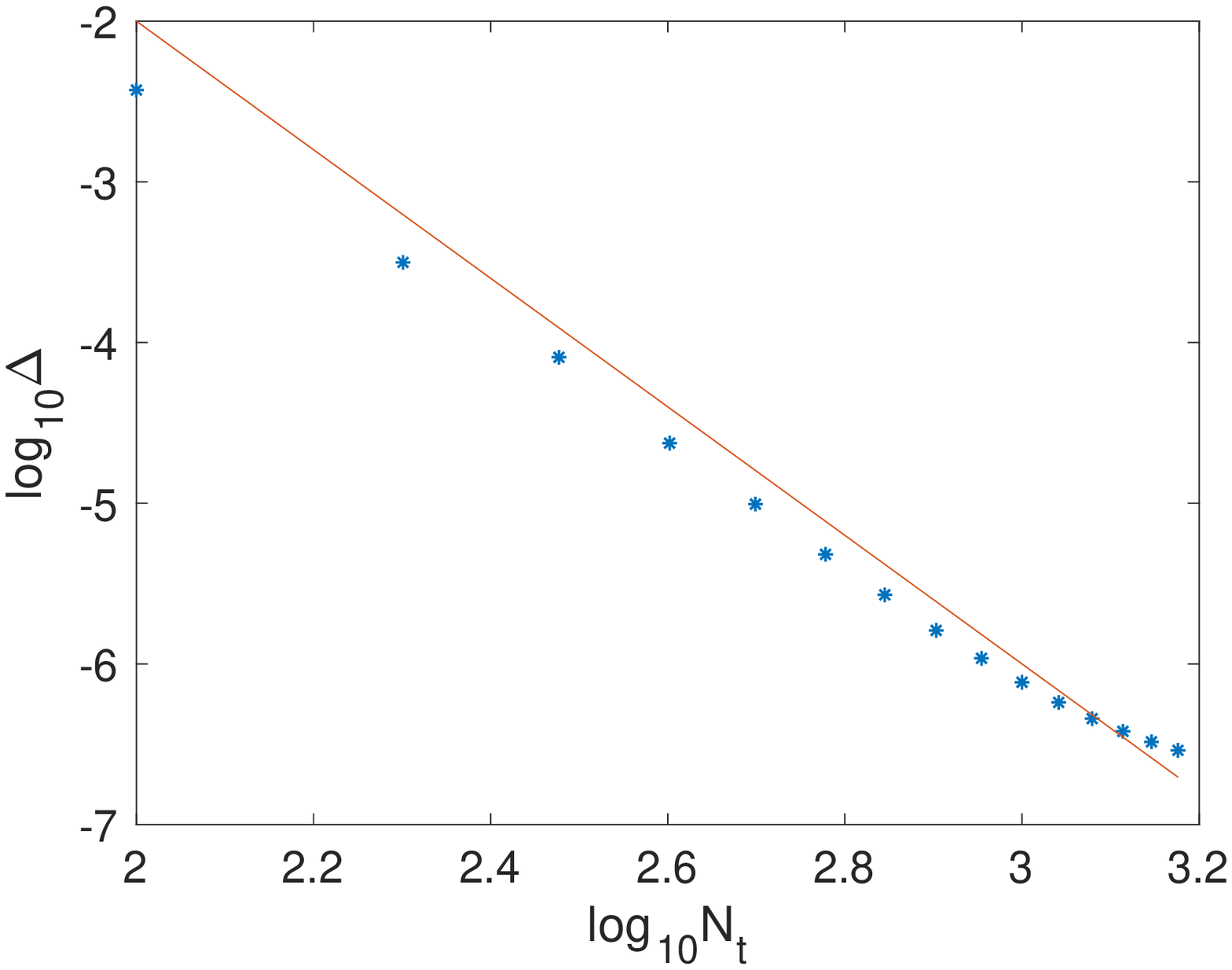}
  \includegraphics[width=0.49\textwidth]{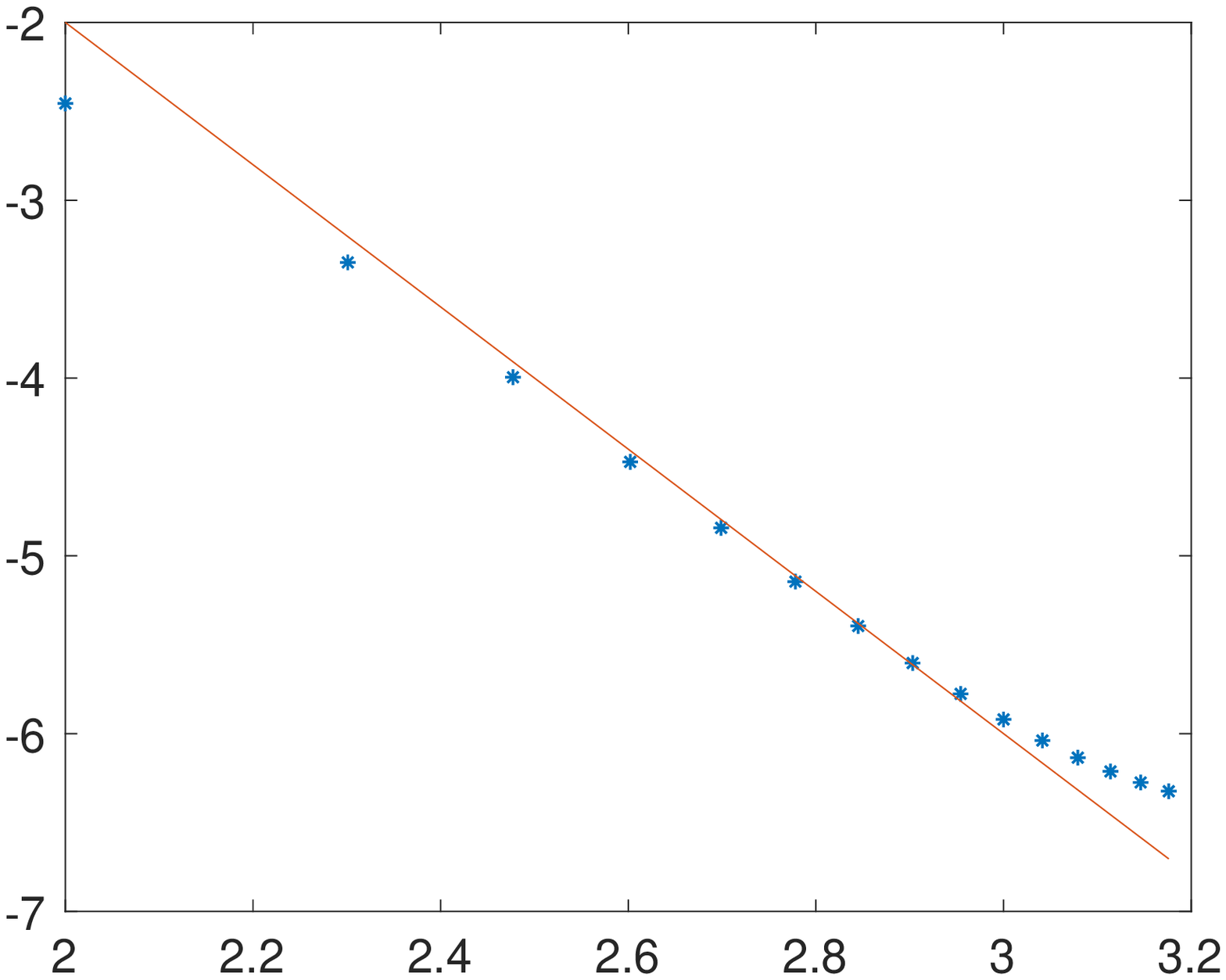}
 \caption{The $L^{\infty}$ norm of the difference 
 between the numerical NLS solution  for Peregrine 
 initial data and the Peregrine solution for $t=1$ on the left, and the 
 numerically computed $E[u]$ (\ref{energy}) on the right; the red line 
 has in both cases slope $-4$.}
 \label{Delta}
\end{figure}

In Fig.~\ref{Delta} we show the $L^{\infty}$ norm of the difference 
between numerical and exact solution for $t=1$ in dependence of the 
number of time steps $N_{t}$. It  can be seen that this is a fourth 
order method by comparing it to a straight line with slope $-4$ in 
the same figure. Below an accuracy of the order of $10^{-6}$ the 
convergence becomes slower than this since the conditioning of the 
Chebyshev differentiation matrices becomes important. This shows that 
numerical errors below $10^{-6}$ can be easily reached, but that for 
higher precisions the code from \cite{BK} based on an iterative 
approach has to be used. On the right of Fig.~\ref{Delta} we show the 
numerically computed quantity (\ref{energy}) in dependence of $N_{t}$. 
It shows a very similar behavior as the error $\Delta$ and is only 
slightly larger. Note that we cannot consider a relative error here 
since $E[u]$ (\ref{energy}) vanishes for the Peregrine solution. 

\section{Localized perturbations}
In this section we consider localized perturbations of the Peregrine 
solution for the elliptic NLS equation ($\kappa=1$). Note that this 
equation does not play a role in the context of the water waves where 
only hyperbolic NLS equations are important, see the discussion in 
\cite{Lannes}. But the elliptic NLS equation plays a role in the 
context of nonlinear optics where the Peregrine solution has been 
recently observed in 1D, see \cite{kibler}. It is also mathematically 
interesting. 

Concretely we study 
the time evolution of the 2D elliptic NLS 
equation for the initial data
\begin{equation}
    u(x,y,0) = u_{Per}(x,0)+c\exp(-(x-x_{0})^{2}-y^{2})
    \label{initial}
\end{equation}
with $c\in\mathbb{C}$ being a constant. We consider the 4 cases 
$x_{0}=0,-1$ 
and $c=\pm 0.1$. All examples are computed with $N_{I}=100$ and 
$N_{II}=101$ polynomials for the $x$-dependence, $N_{y}=2^{7}$ 
Fourier modes for 
$y=3[-\pi,\pi]$ and $N_{t}=1000$ time steps. Note that larger values 
of these parameters have only an effect on the solution below 
plotting accuracy, i.e., the resulting figures would be 
indistinguishable from the ones shown. 

\begin{figure}[htb!]
 \includegraphics[width=0.49\textwidth]{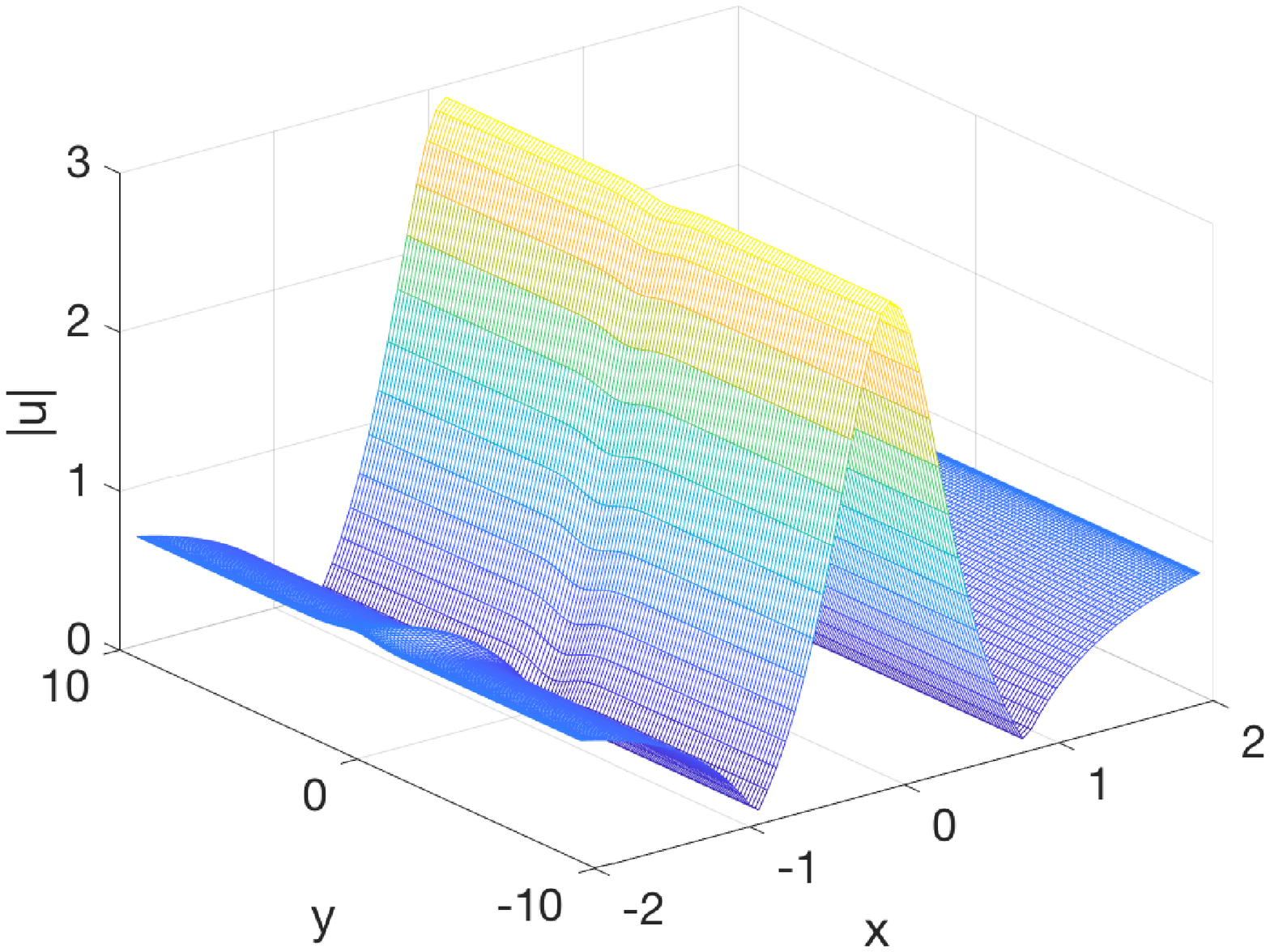}
   \includegraphics[width=0.49\textwidth]{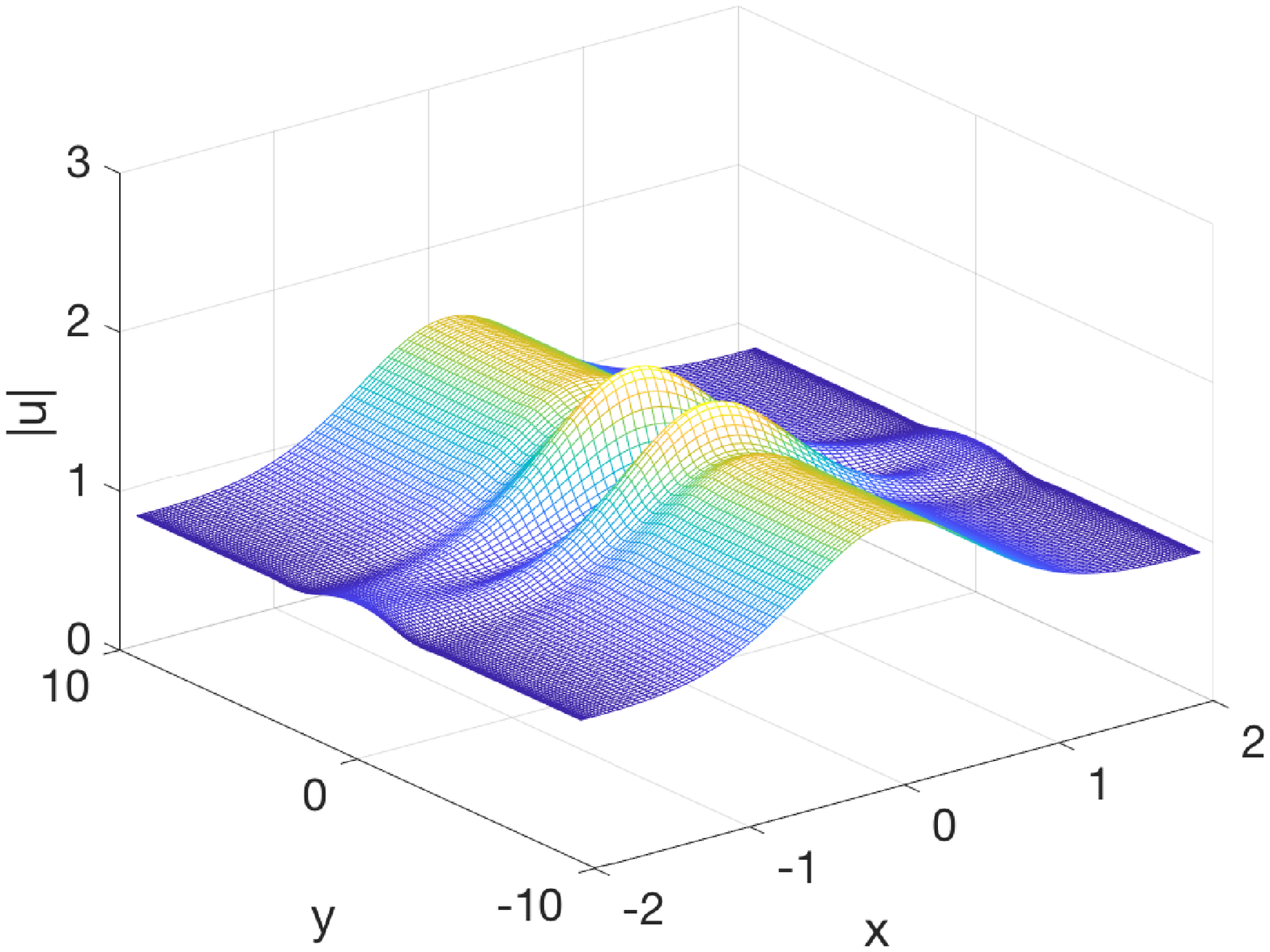}
 \caption{Solution to the 2D NLS equation for the initial data 
 $u(x,y,0)=u_{Per}(x,t_{0})+0.1\exp(-(x+1)^{2}-y^{2})$ for $t=0$ on the 
 left and $t=0.5$ on the right.}
 \label{Peregine2dp01gaussx1}
\end{figure}

We first consider the case $x_{0}=-1$ and $c=0.1$ in 
Fig.~\ref{Peregine2dp01gaussx1}, on the left the initial condition, 
on the right the solution for $t=0.5$. The solution is clearly 
unstable in the sense that the initial perturbations grow.

\begin{figure}[htb!]
  \includegraphics[width=0.49\textwidth]{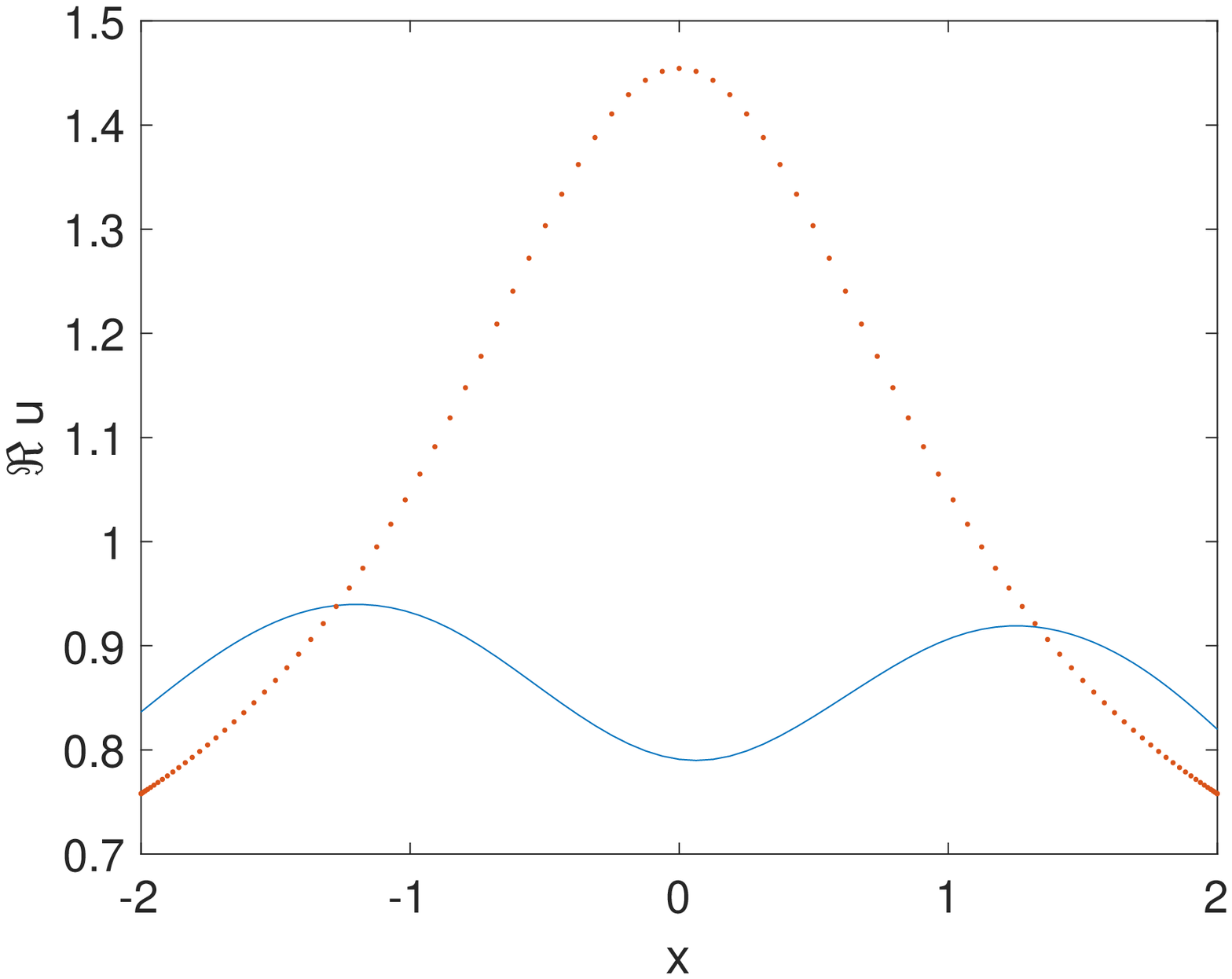}
  \includegraphics[width=0.49\textwidth]{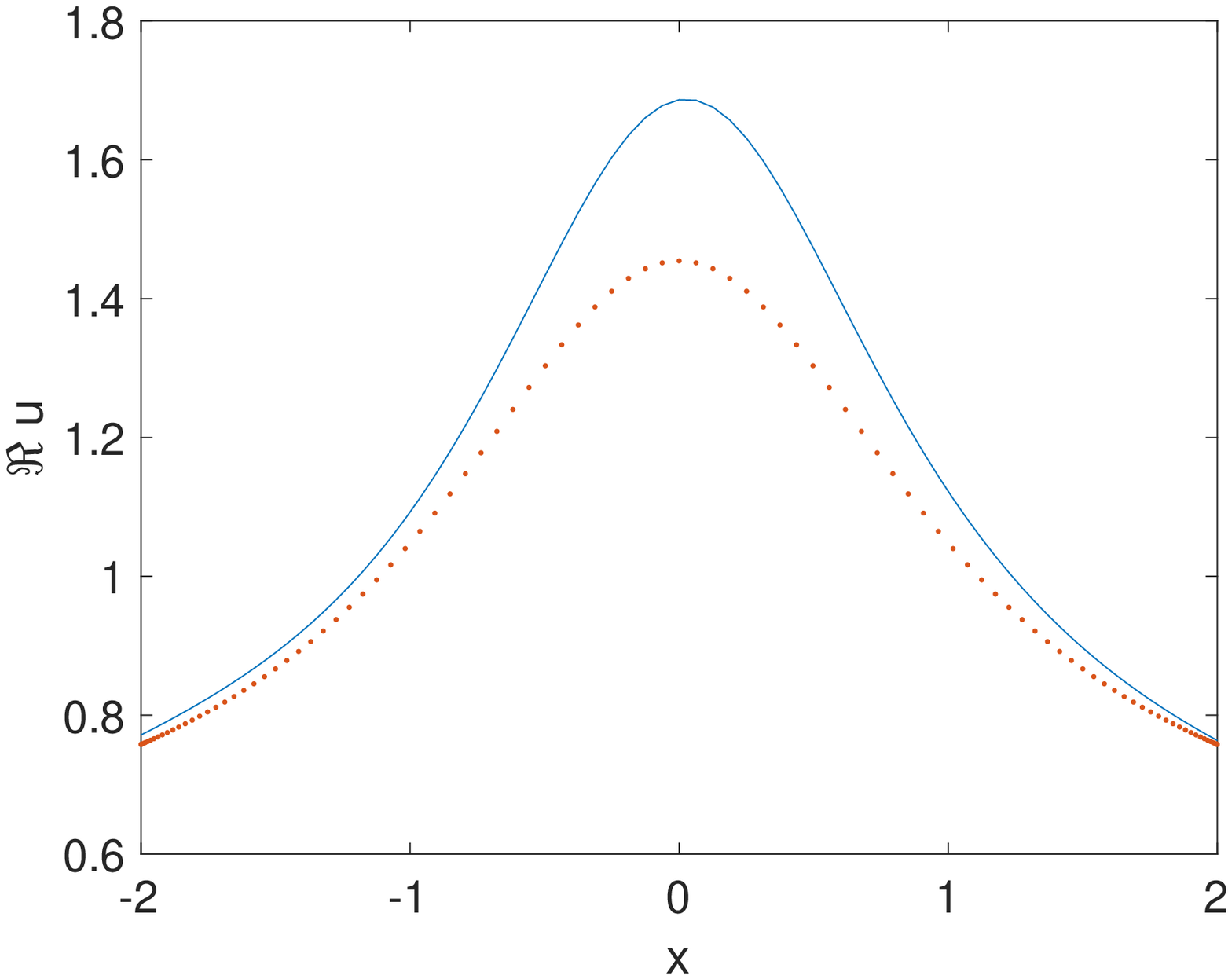}
 \caption{Real part of the solution to the 2D NLS equation for the initial data 
 $u(x,y,0)=u_{Per}(x,t_{0})+0.1\exp(-(x+1)^{2}-y^{2})$ for $t=0.5$,  on the 
 left for $y=0$, on the right for $y=1.6199$; the Peregrine solution 
 for the same time is shown as a dotted line.}
 \label{Peregine2dp01gaussx1y}
\end{figure}

To show that the solution does not stay close to the Peregrine 
solution, we show in Fig.~\ref{Peregine2dp01gaussx1y} the solutions 
for $t=0.5$ and for $y=0$  (minimal modulus of $u$ for $x=0$ and 
$t=0.5$) and 
$y=1.6199$ (one of the maxima of $|u|$ for $x=0$ and $t=0.5$) 
together with the Peregrine solution. 

\begin{figure}[htb!]
  \includegraphics[width=0.49\textwidth]{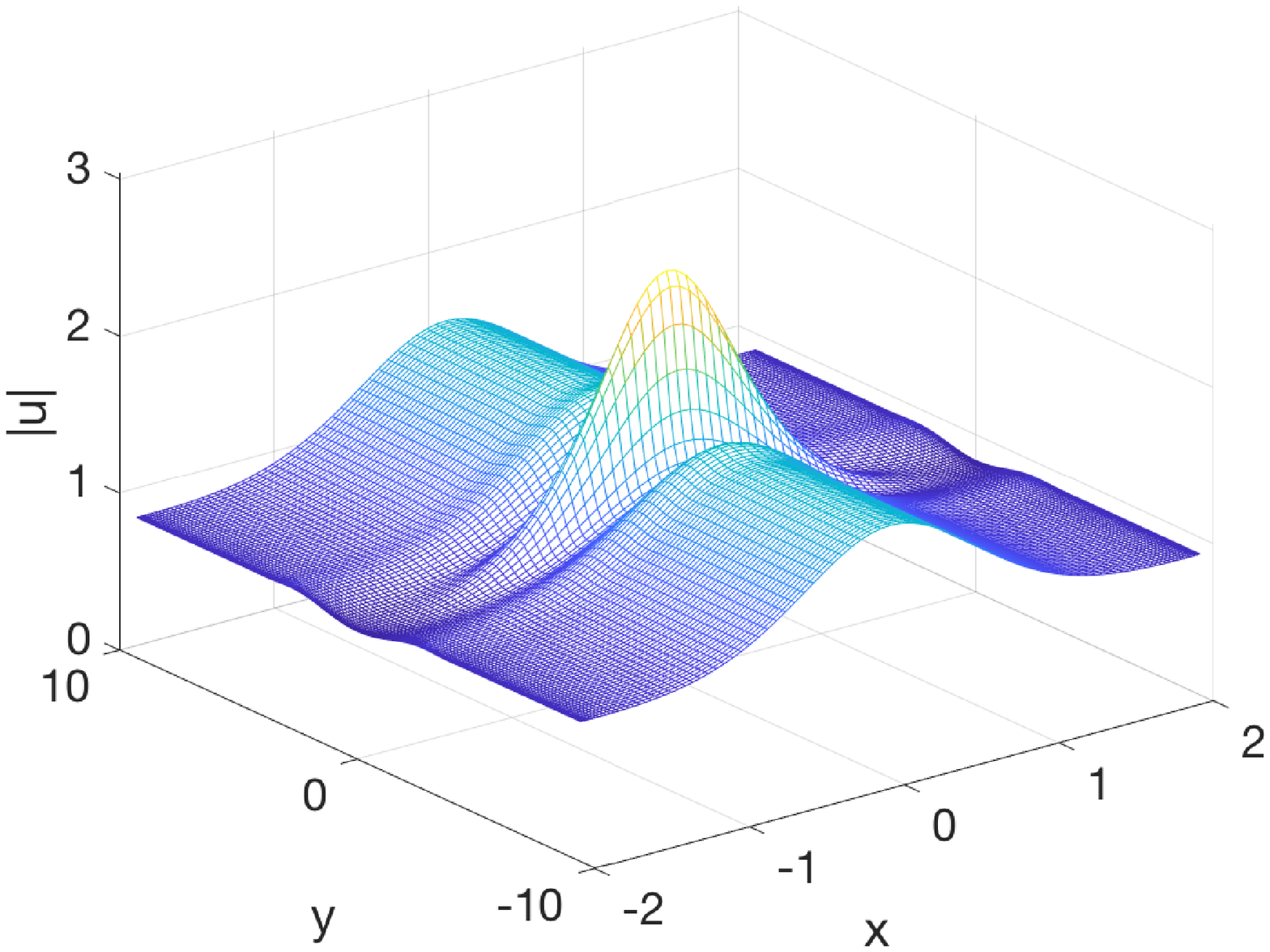}
  \includegraphics[width=0.49\textwidth]{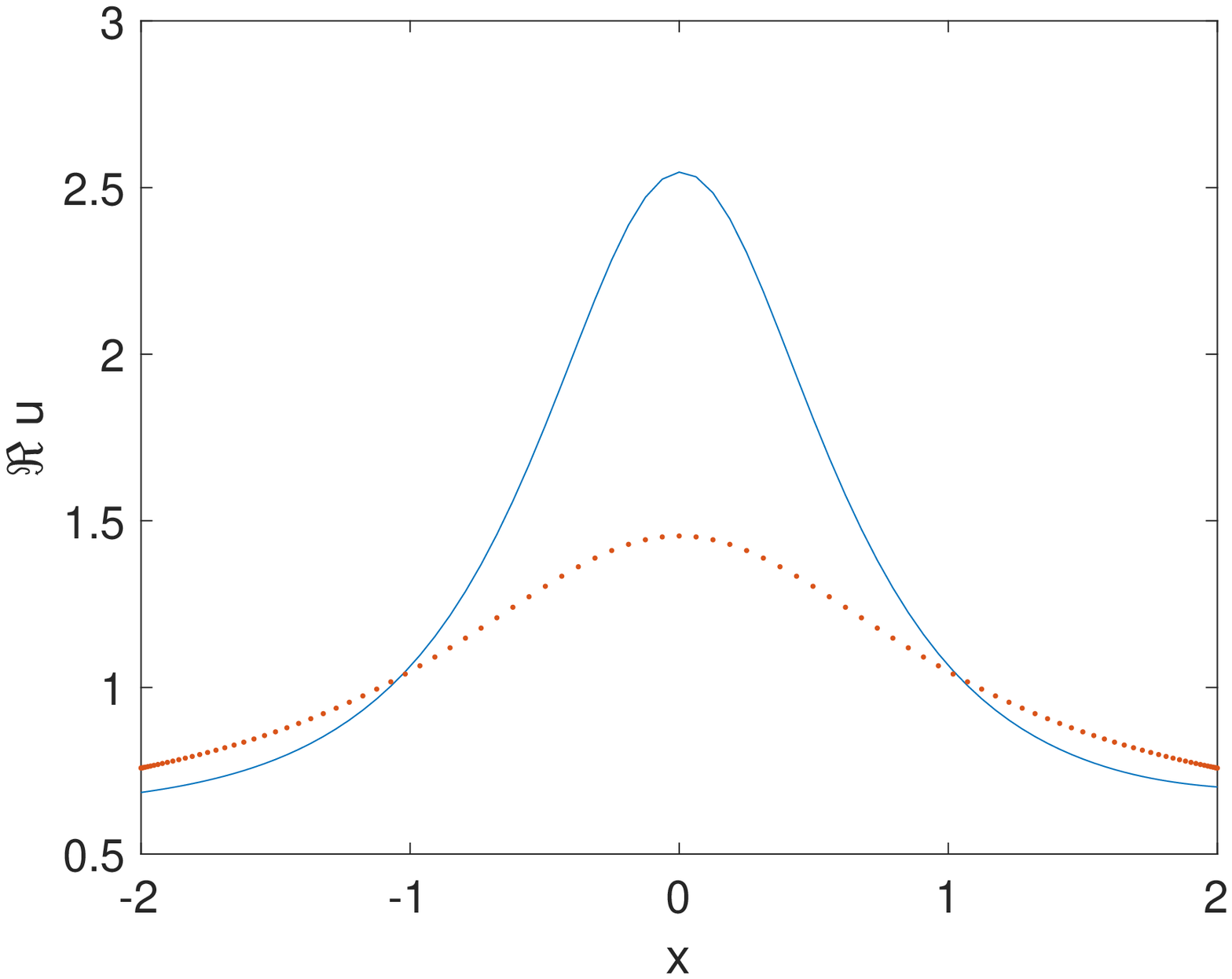}
 \caption{Solution to the 2D NLS equation for the initial data 
 $u(x,y,0)=u_{Per}(x,t_{0})-0.1\exp(-(x+1)^{2}-y^{2})$ for $t=0.5$ on the 
 left, and the real part of the solution at the same time for $y=0$ on the right together with the 
 Peregrine solution (dotted).}
 \label{Peregine2dm01gaussx1}
\end{figure}

The situation is similar in the case $x_{0}=-1$ and $c=-0.1$ as can be 
seen in Fig.~\ref{Peregine2dm01gaussx1}: on the left the modulus of 
the solution for $t=0.5$ is shown, on the right the solution for 
$y=0$ together with the Peregrine solution (dotted). Again the 
perturbed solution does not stay close to the Peregrine solution. 

\begin{figure}[htb!]
 \includegraphics[width=0.49\textwidth]{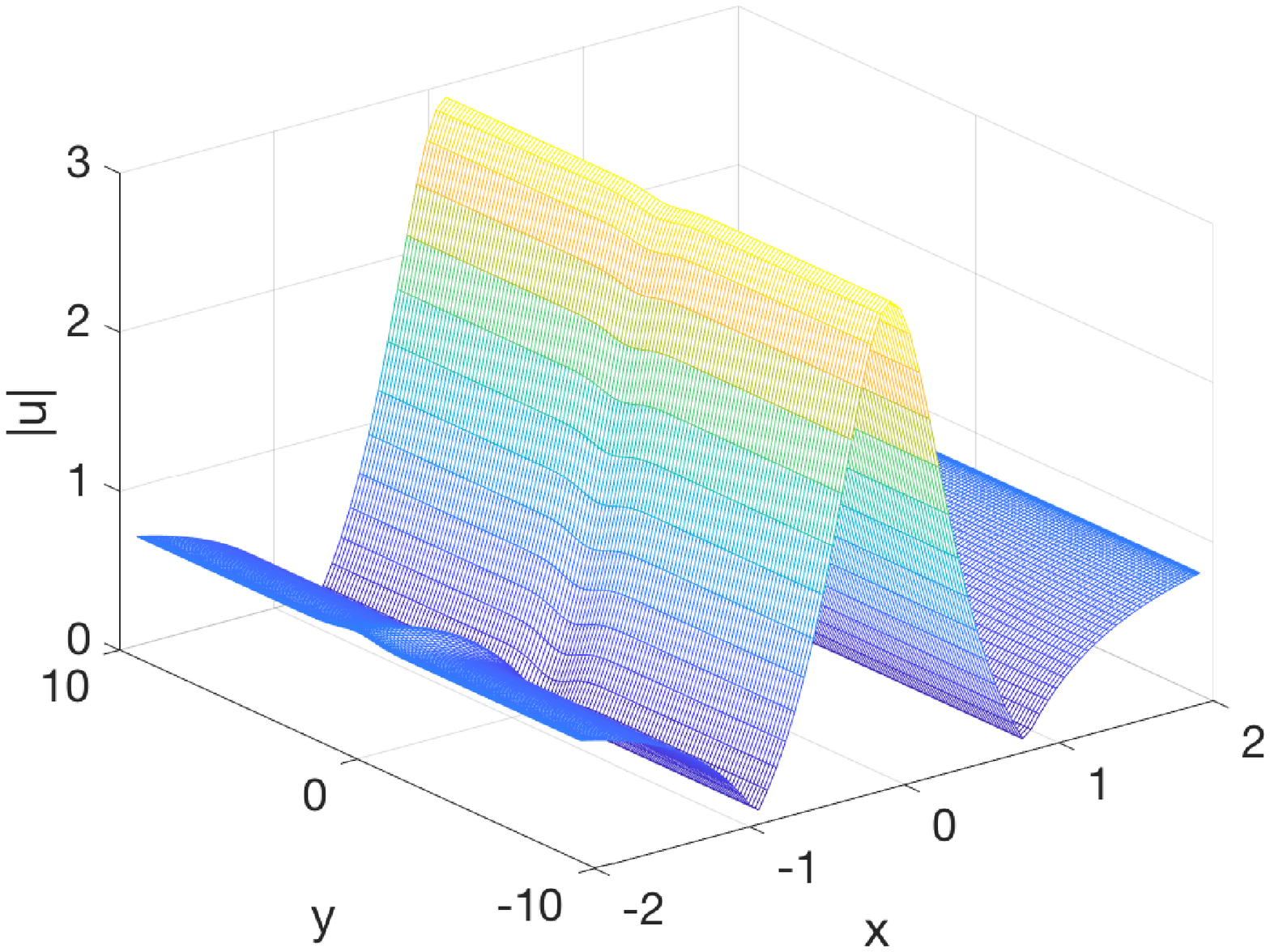}
   \includegraphics[width=0.49\textwidth]{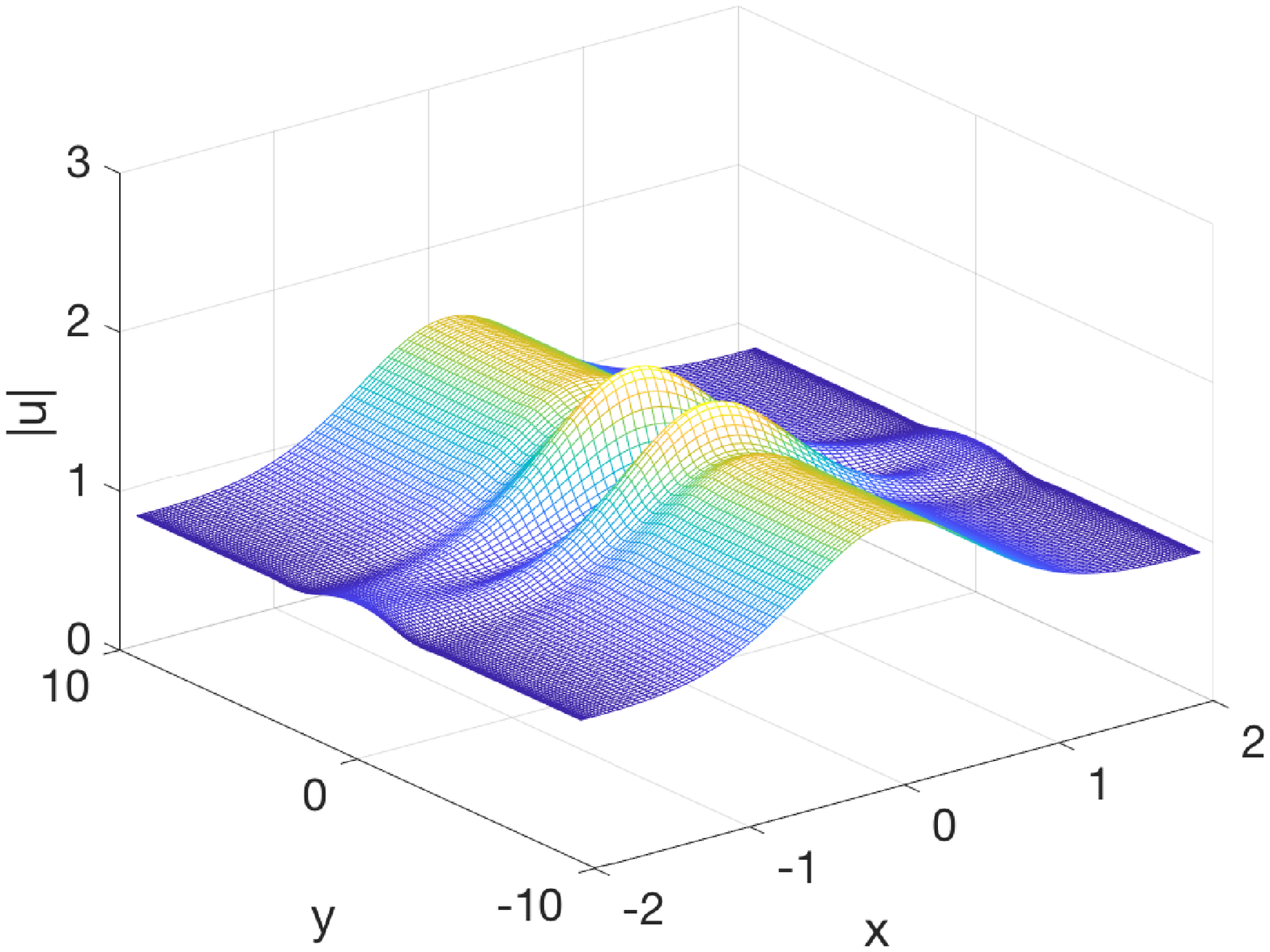}
 \caption{Solution to the 2D NLS equation for the initial data 
 $u(x,y,0)=u_{Per}(x,t_{0})-0.1\exp(-(x-1)^{2}-y^{2})$ for $t=0$ on the 
 left and $t=0.5$ on the right.}
 \label{Peregine2dp01gauss}
\end{figure}

The solution for the initial data (\ref{initial}) for $x_{0}=0$ and 
$c=0.1$ (on the left of Fig.~\ref{Peregine2dp01gauss}) at time $t=0.5$ 
can be seen on the right of Fig.~\ref{Peregine2dp01gauss}. The 
initial perturbations appear once more to grow in time. 

\begin{figure}[htb!]
  \includegraphics[width=0.49\textwidth]{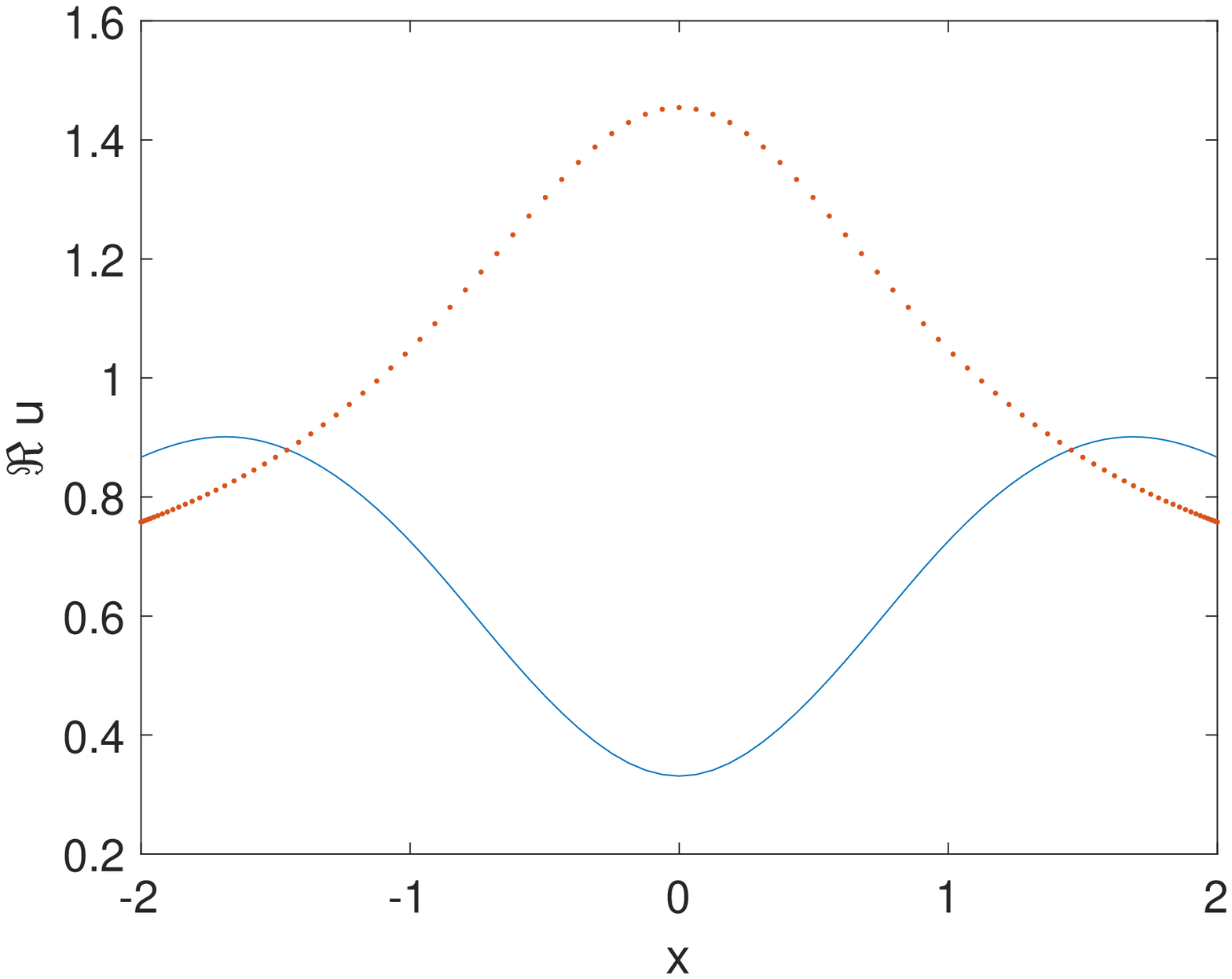}
  \includegraphics[width=0.49\textwidth]{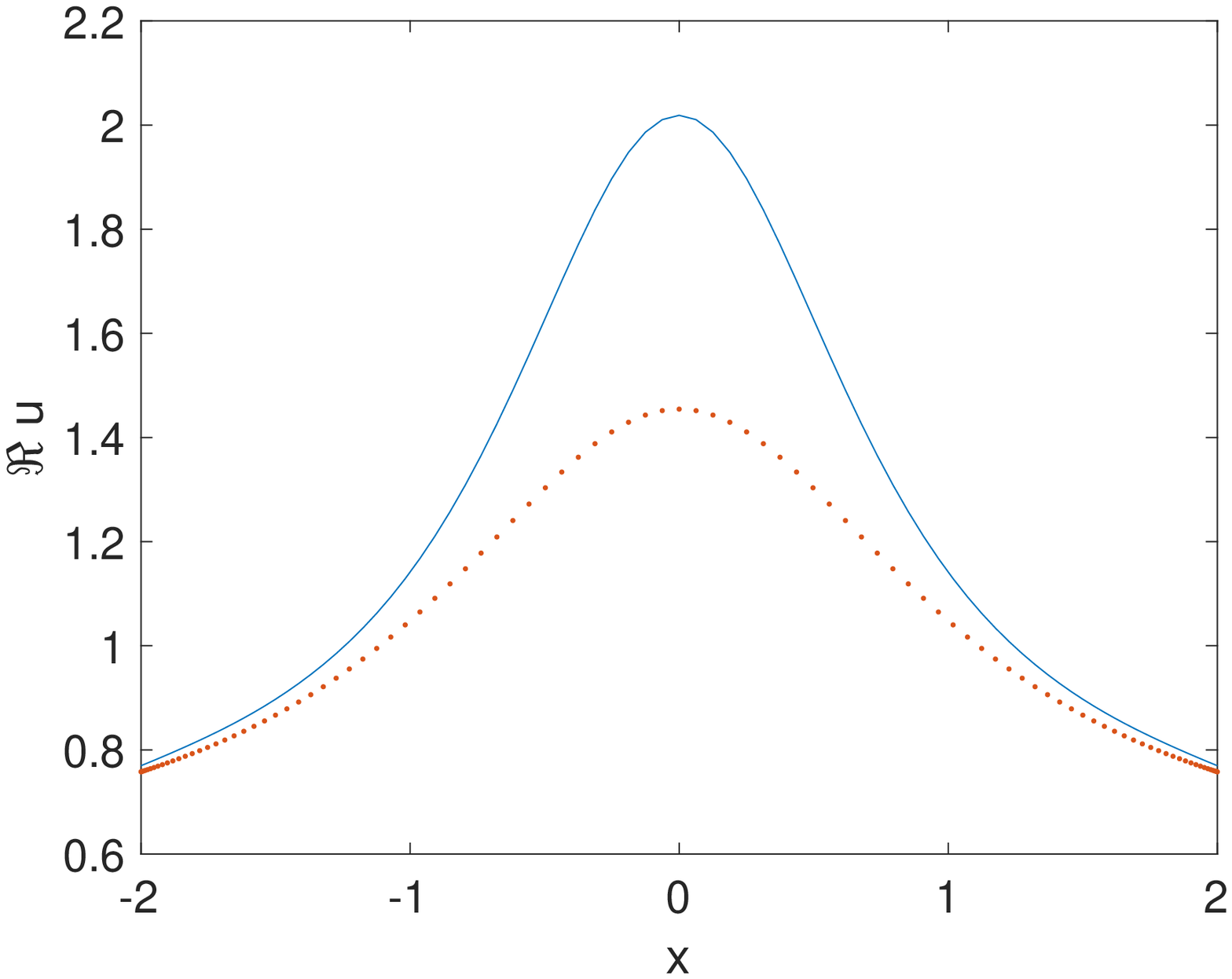}
 \caption{Real part of the solution to the 2D NLS equation for the initial data 
 $u(x,y,0)=u_{Per}(x,t_{0})+0.1\exp(-x^{2}-y^{2})$ for $t=0.5$, on the 
 left for $y=0$ and on the right for  $y=1.7671$; the Peregrine solution 
 for the same time is shown as a dotted line.}
 \label{Peregine2dp01gaussy}
\end{figure}

In Fig.~\ref{Peregine2dp01gaussy} we show the solution on the right 
of Fig.~\ref{Peregine2dp01gauss} for $y=0$ and for $y=1.7671$, one of 
the maxima of the modulus of $u$ for $x=0$ together with the Peregrine 
solution. 

\begin{figure}[htb!]
  \includegraphics[width=0.49\textwidth]{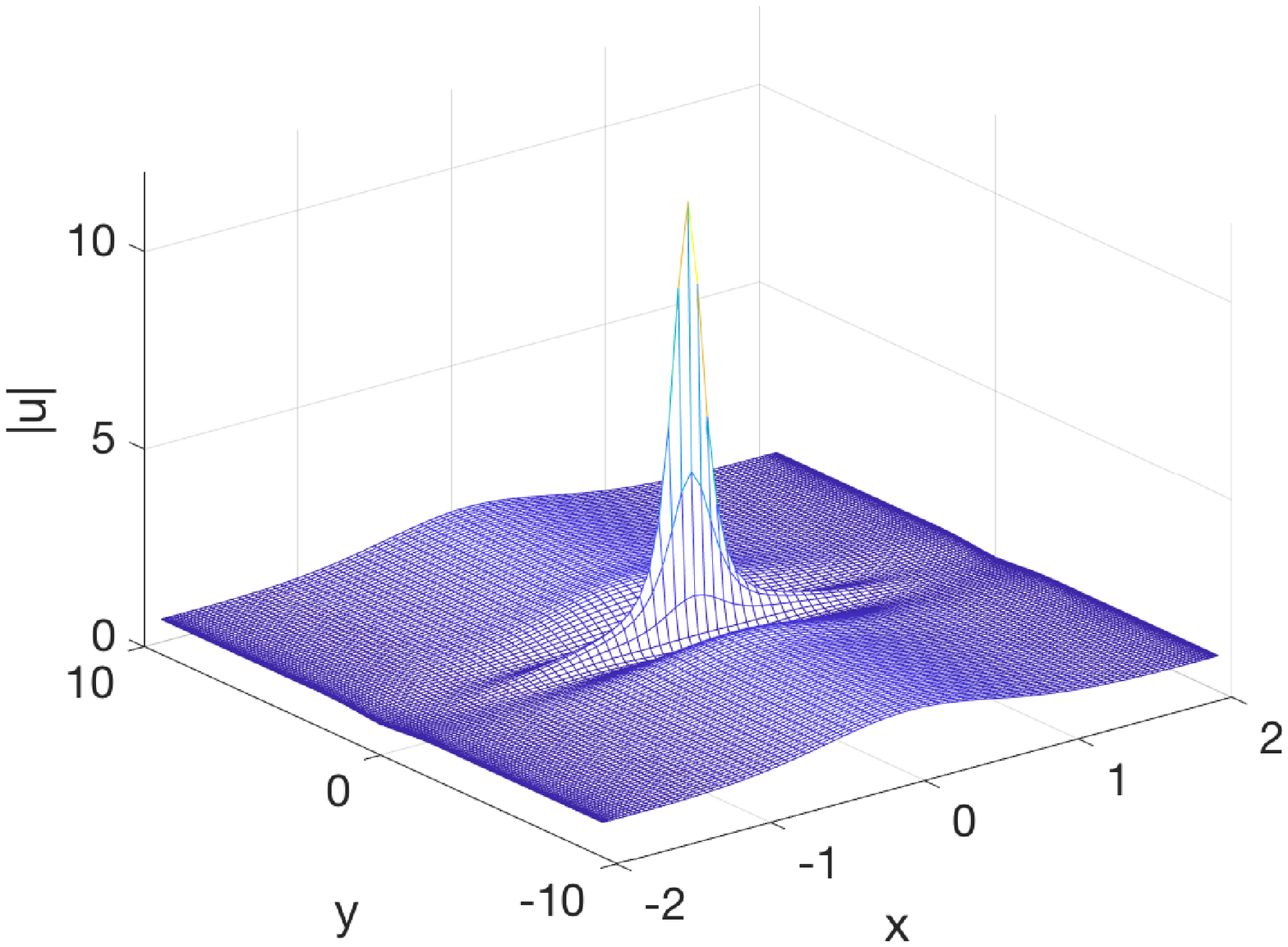}
  \includegraphics[width=0.49\textwidth]{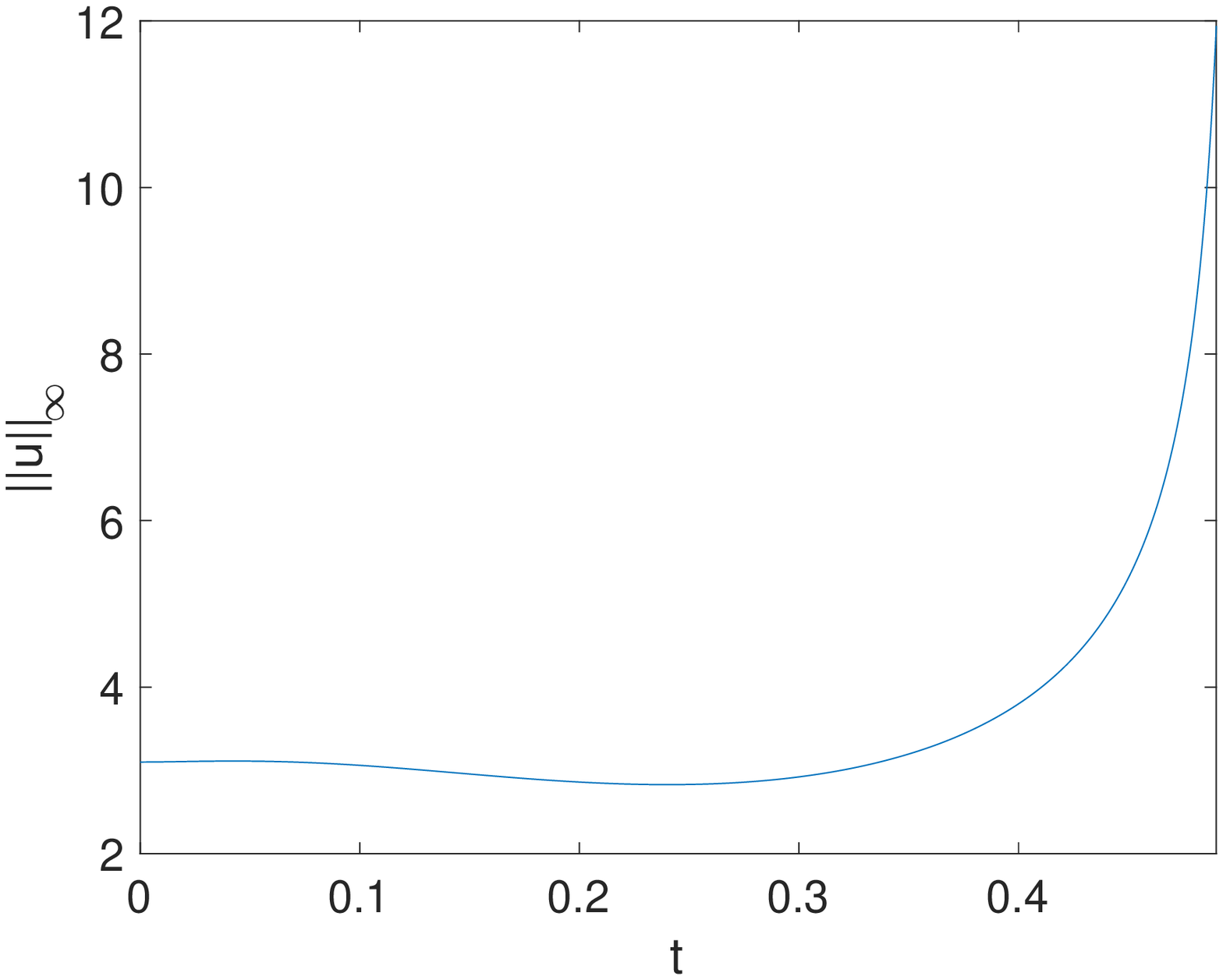}
 \caption{Solution to the 2D NLS equation for the initial data 
 $u(x,y,0)=u_{Per}(x,t_{0})-0.1\exp(-x^{2}-y^{2})$ for $t=0.49$ on the 
 left and the $L^{\infty}$ norm of the solution in dependence of 
 time on the right.}
 \label{Peregine2dm01gauss}
\end{figure}

The situation changes somewhat in the case $x_{0}=0$ and $c=-0.1$. 
In this case the solution appears to blow up for $t\sim0.49$ (the 
relative $E[u]$ conservation drops in this case below $10^{-3}$ and 
the code is stopped). The solution for $t=0.49$ can be seen on the 
left of Fig.~\ref{Peregine2dm01gauss}. The $L^{\infty}$ norm 
of the solution in dependence of time on the 
right of the same figure  also indicates a blow-up. 
This would imply that the Peregrine solution is strongly unstable as 
a solution to the 2D focusing elliptic NLS equation. 

\section{Nonlocalized perturbations}
In this section we study perturbations  nonlocalized in $x$, but 
localized in $y$ of the 
Peregrine solution $u_{Per}$ of the form $u(x,t_{0})=(\sigma-\exp(-y^{2})) u_{Per}(x,t_{0})$, where 
$\sigma\sim 1$. The same numerical parameters as in the previous 
section are applied. We concentrate again on the elliptic case 
$\kappa=1$. 

\begin{figure}[htb!]
  \includegraphics[width=0.49\textwidth]{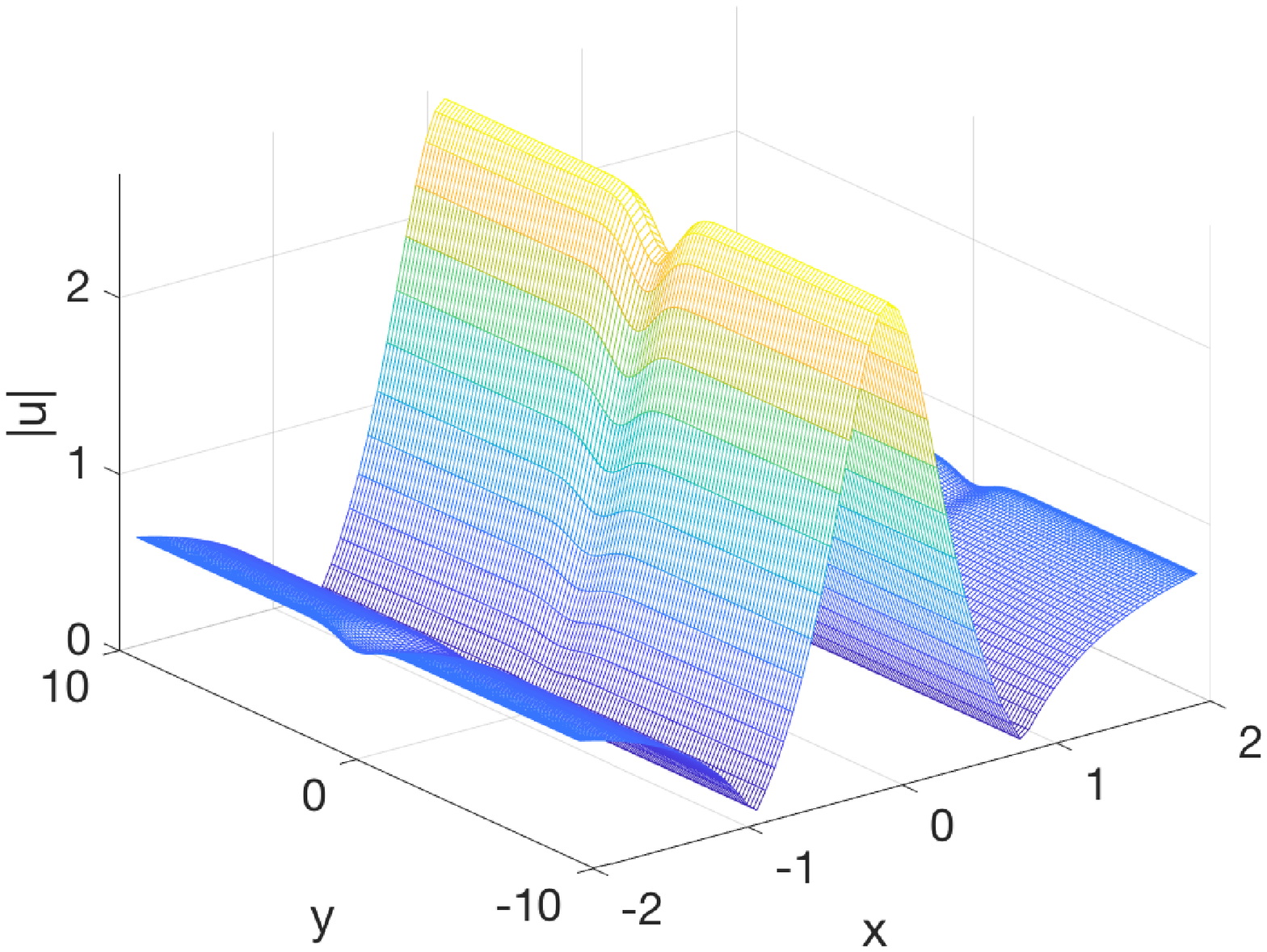}
  \includegraphics[width=0.49\textwidth]{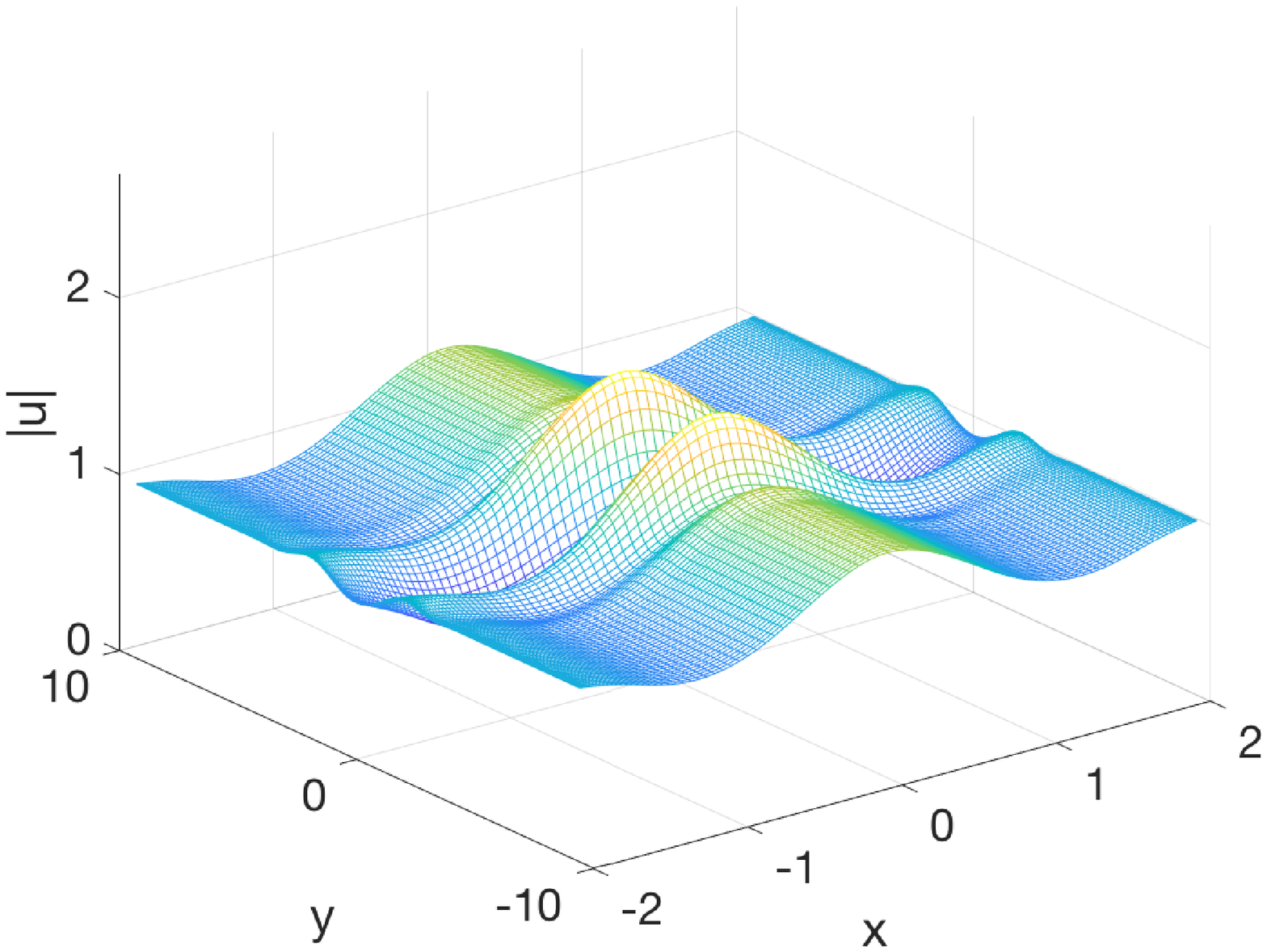}
 \caption{Solution to the 2D NLS equation for the initial data 
 $(0.9-\exp(-y^{2})) u_{Per}(x,t_{0})$ for $t=0$ on the left and for 
 $t=0.5$ on the right.}
 \label{Peregine2d09}
\end{figure}

The solution for the case $\sigma=0.9$ (in Fig.~\ref{Peregine2d09} on 
the left) can be seen for $t=0.5$ in 
Fig.~\ref{Peregine2d09} on the right. Clearly the initial 
perturbations grow.

\begin{figure}[htb!]
  \includegraphics[width=0.49\textwidth]{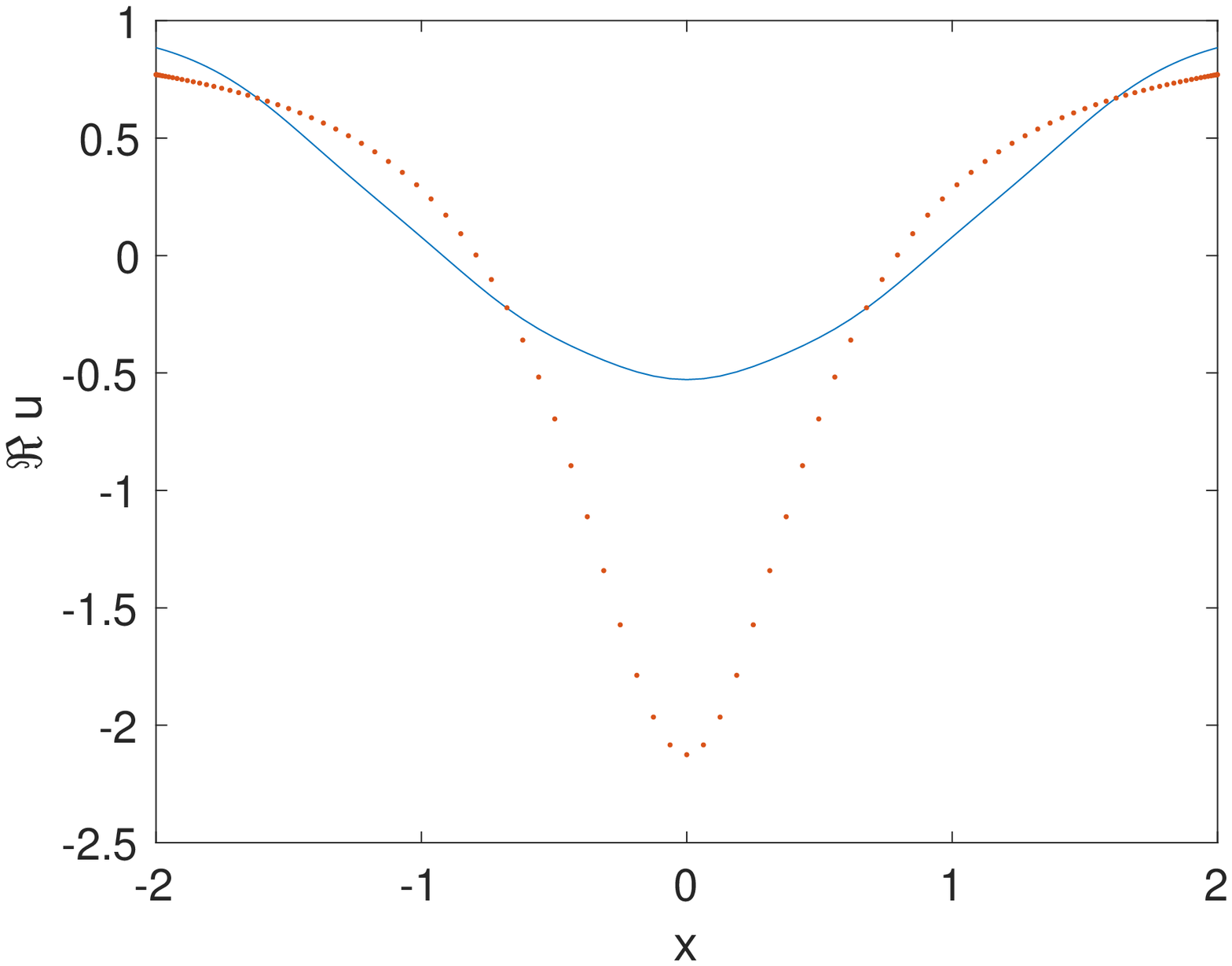}
  \includegraphics[width=0.49\textwidth]{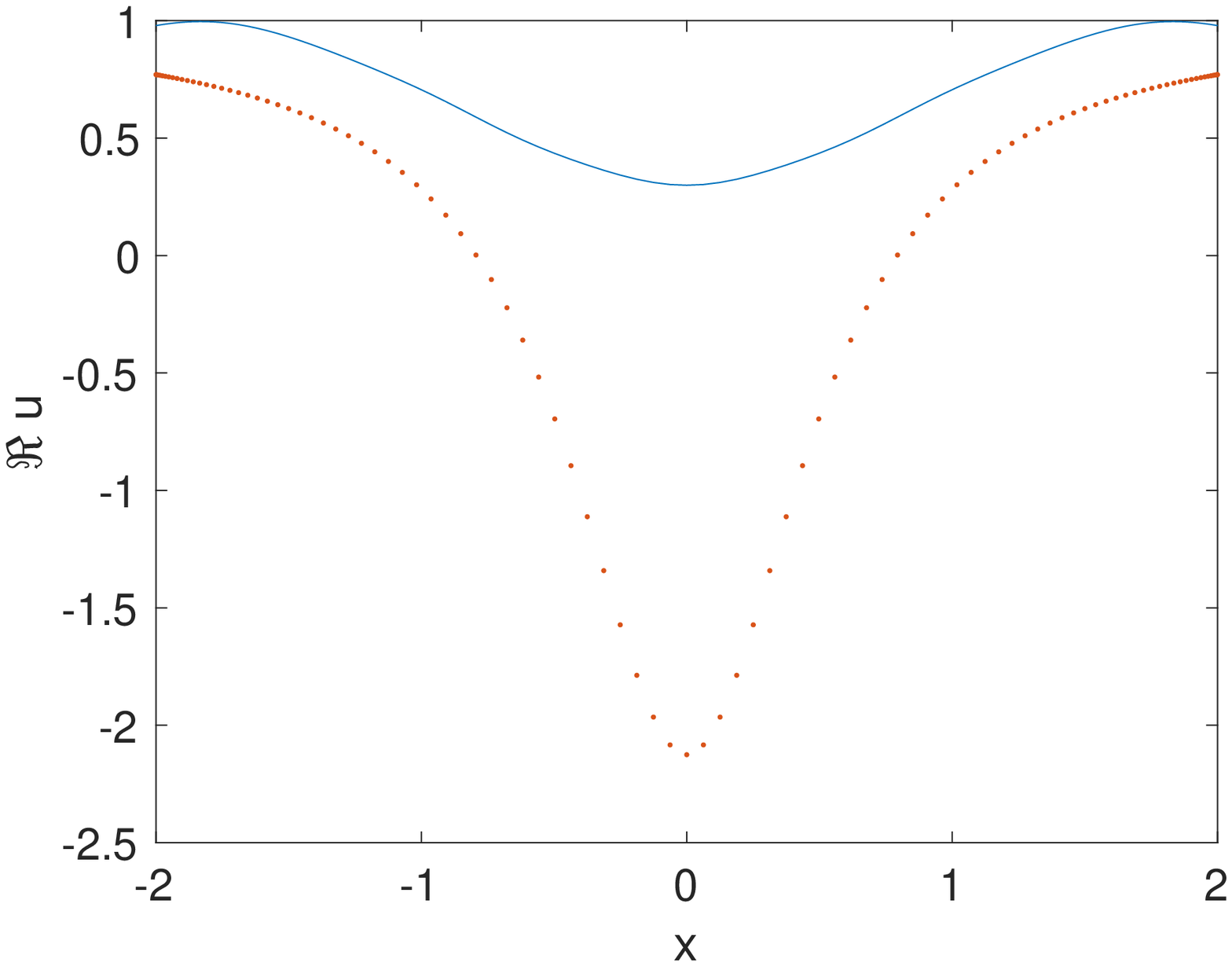}
 \caption{Real part of the solution to the 2D NLS equation for the initial data 
 $(0.9-\exp(-y^{2})) u_{Per}(x,t_{0})$ for $t=0.5$ together with the 
 dotted Peregrine solution, on the left for 
 $y=0$, on the tight for $y=2.0617$.}
 \label{Peregine2d09y}
\end{figure}

In Fig.~\ref{Peregine2d09y} we show the solution on the right of 
Fig.~\ref{Peregine2d09} for $y=0$ and $y=2.0617$ (the maximum of the 
modulus of $|u|$) together with the 
dotted Peregrine solution. It can be seen that the perturbed NLS solution does not 
stay close to the latter. 

\begin{figure}[htb!]
  \includegraphics[width=0.49\textwidth]{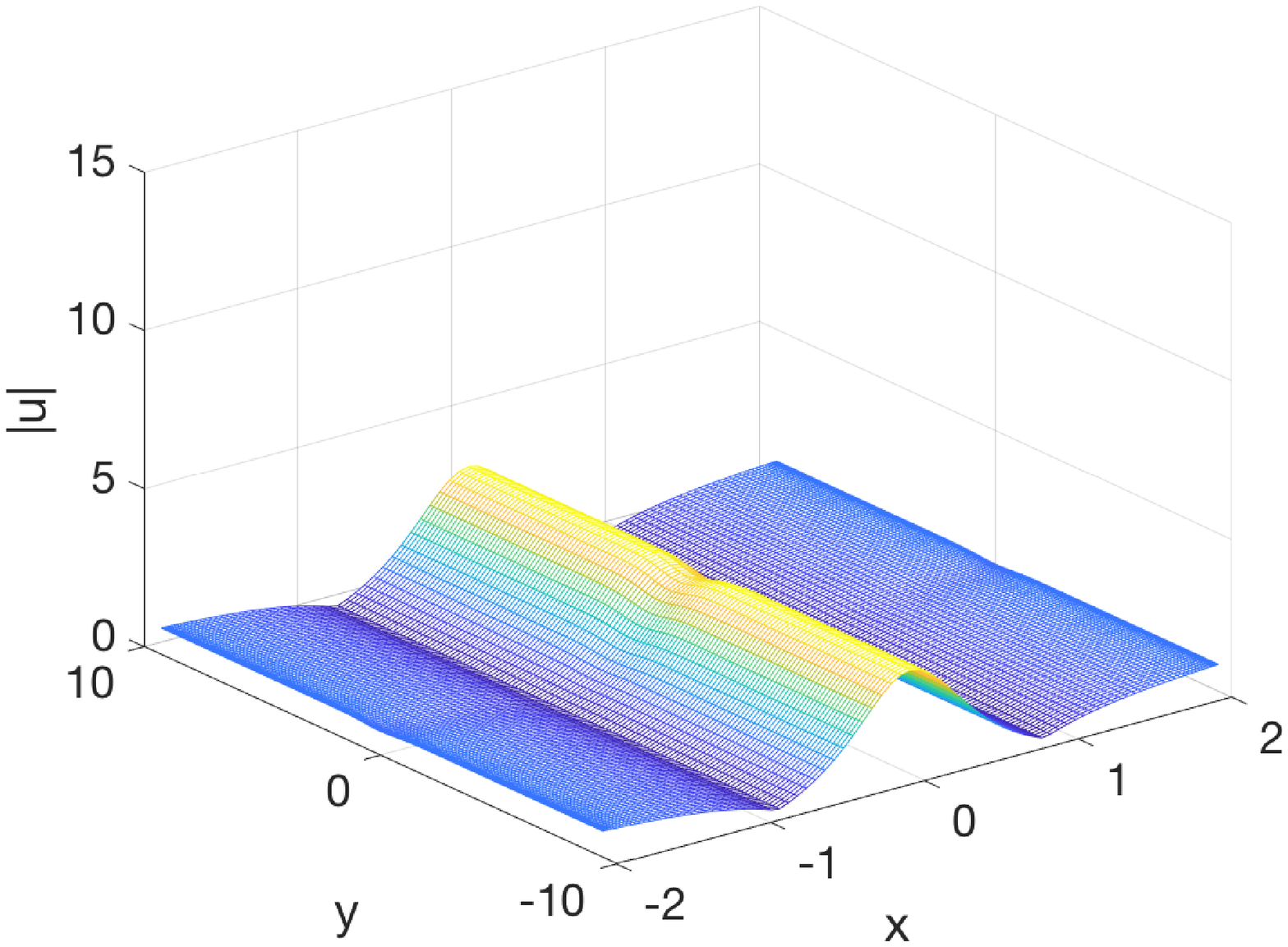}
  \includegraphics[width=0.49\textwidth]{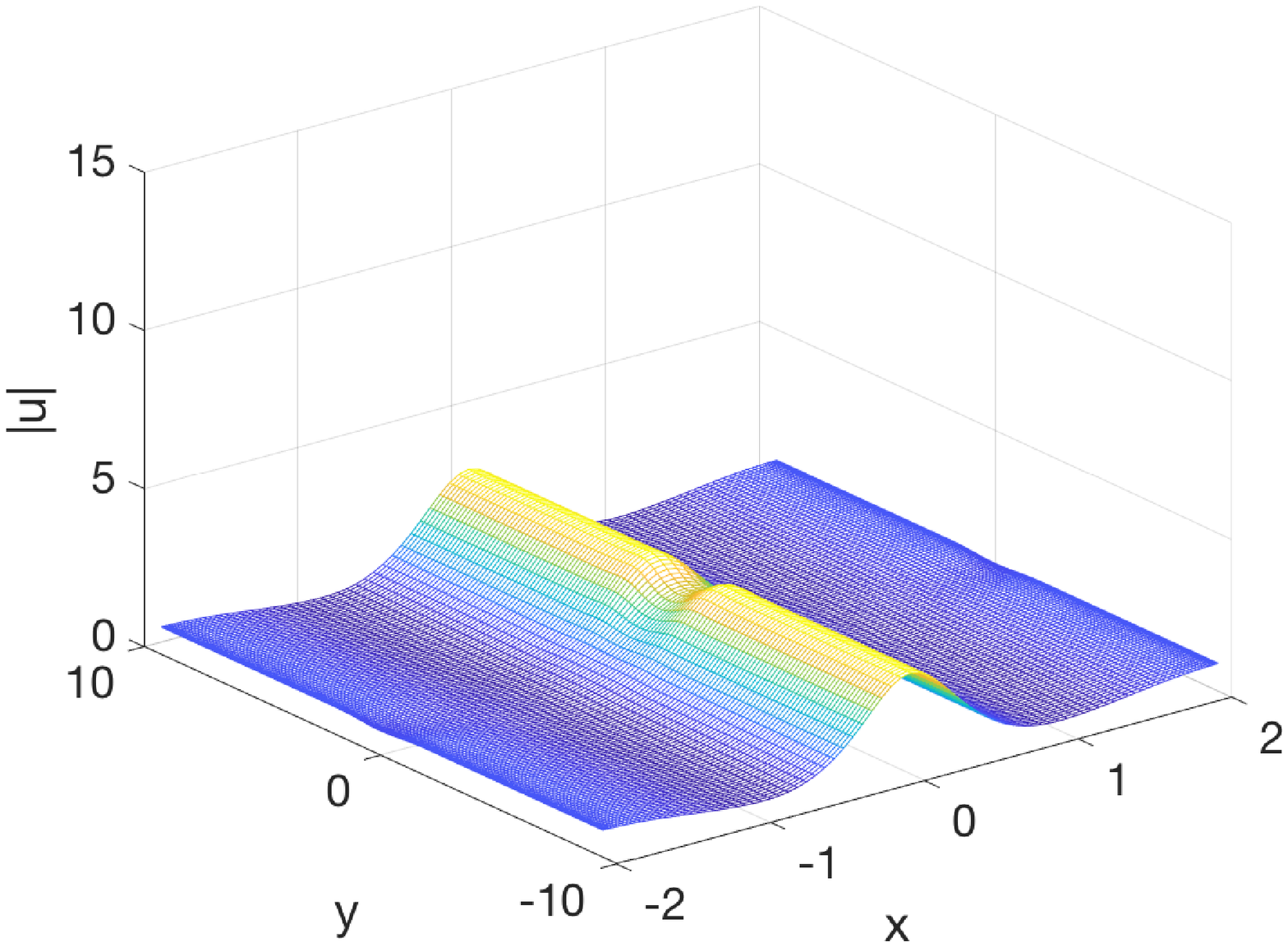}\\
    \includegraphics[width=0.49\textwidth]{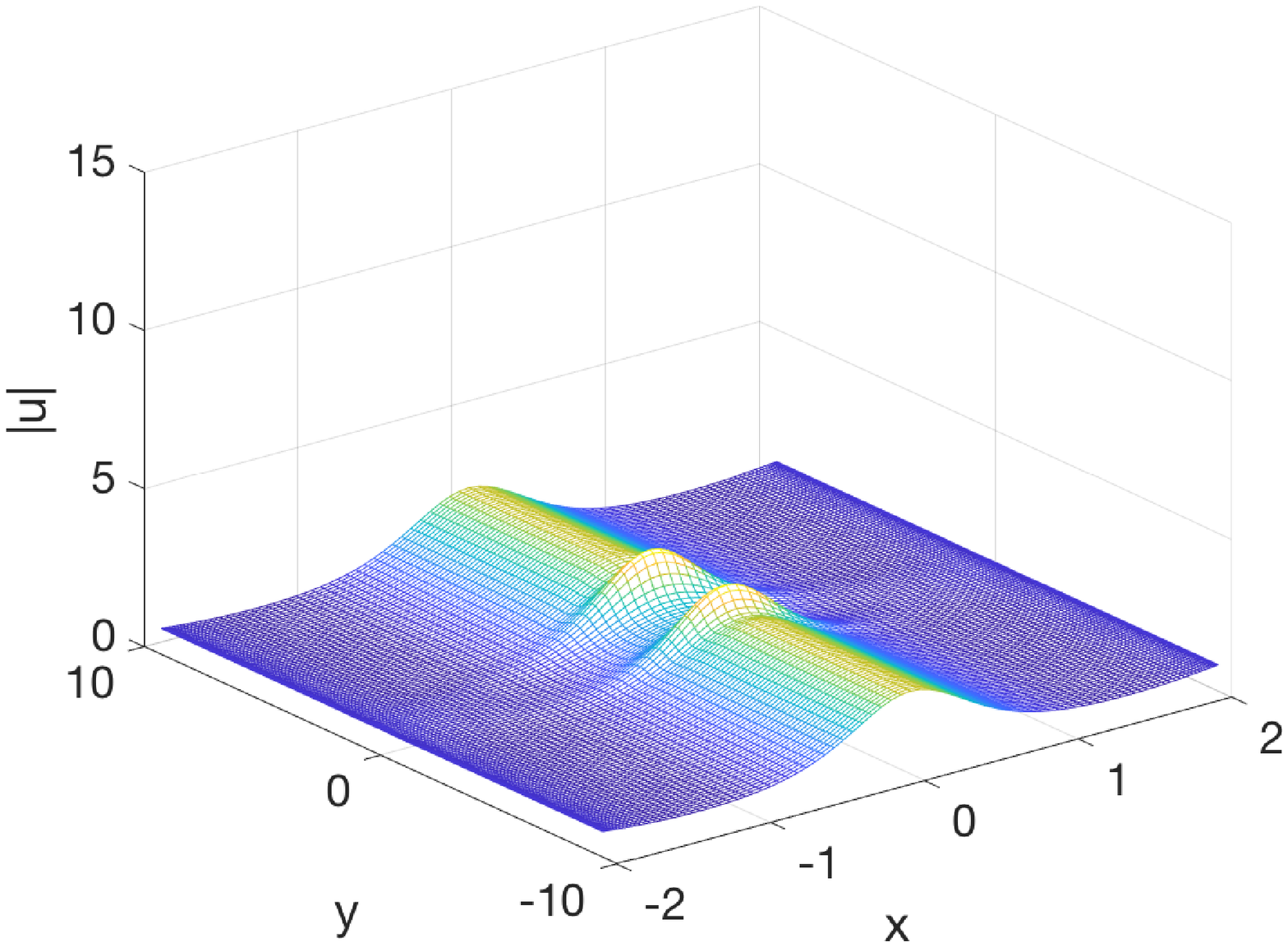}
  \includegraphics[width=0.49\textwidth]{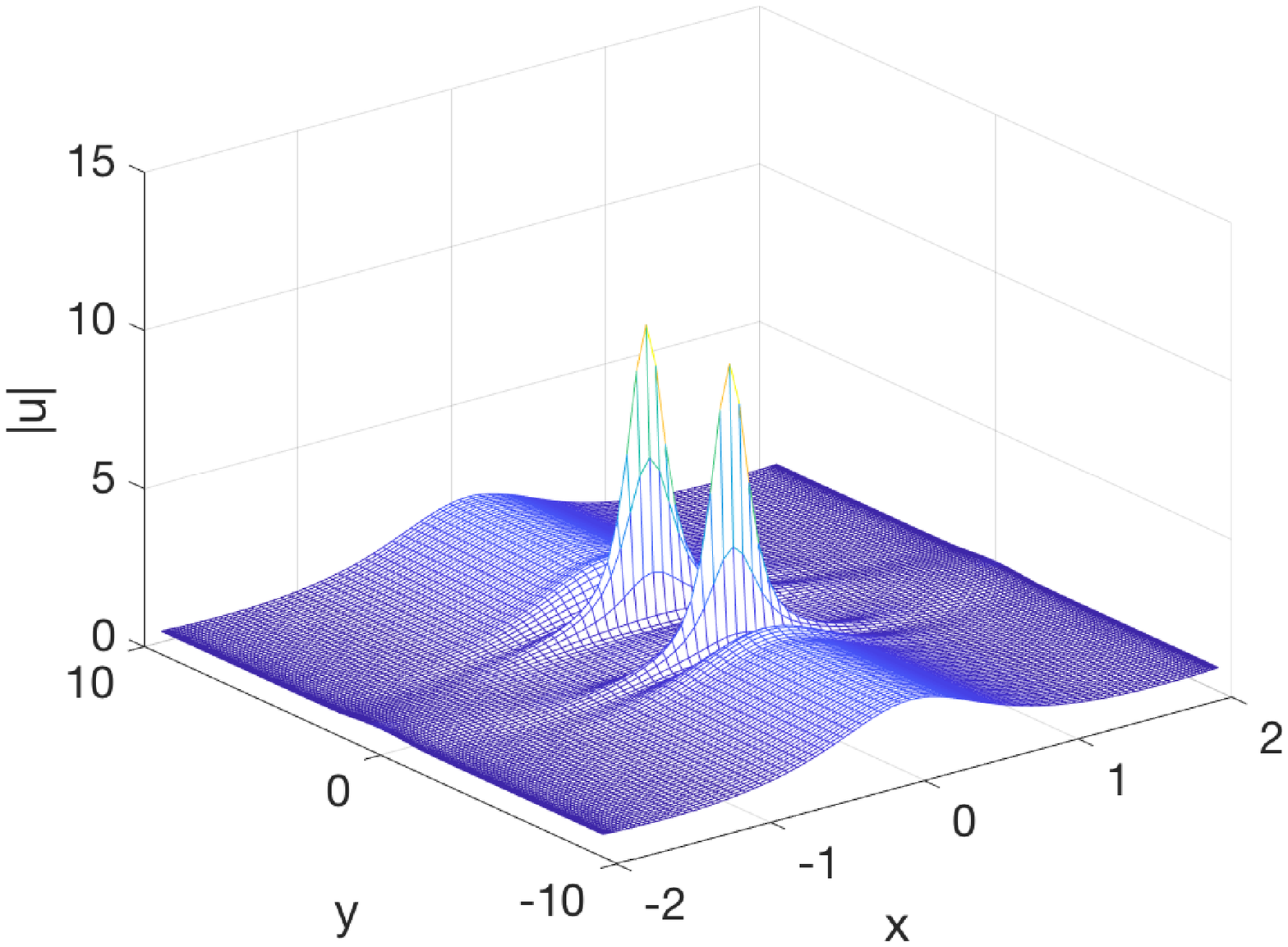}
 \caption{Solution to the 2D NLS equation for the initial data 
 $(1.1-\exp(-y^{2})) u_{Per}(x,t_{0})$ for $t=0$, $t=0.134$, (left 
 respectively right in the upper 
 row), and 
 $t=0.268$ and $t=0.4$ (left 
 respectively right in the lower row).}
 \label{Peregine2d11}
\end{figure}

The solution for the case $\sigma=1.1$ can be seen at different times 
in Fig.~\ref{Peregine2d11}. The initial 
perturbation obviously grows, and the formed structures appear even 
to blow up in this case for $t\sim 0.4$. 

\begin{figure}[htb!]
  \includegraphics[width=0.49\textwidth]{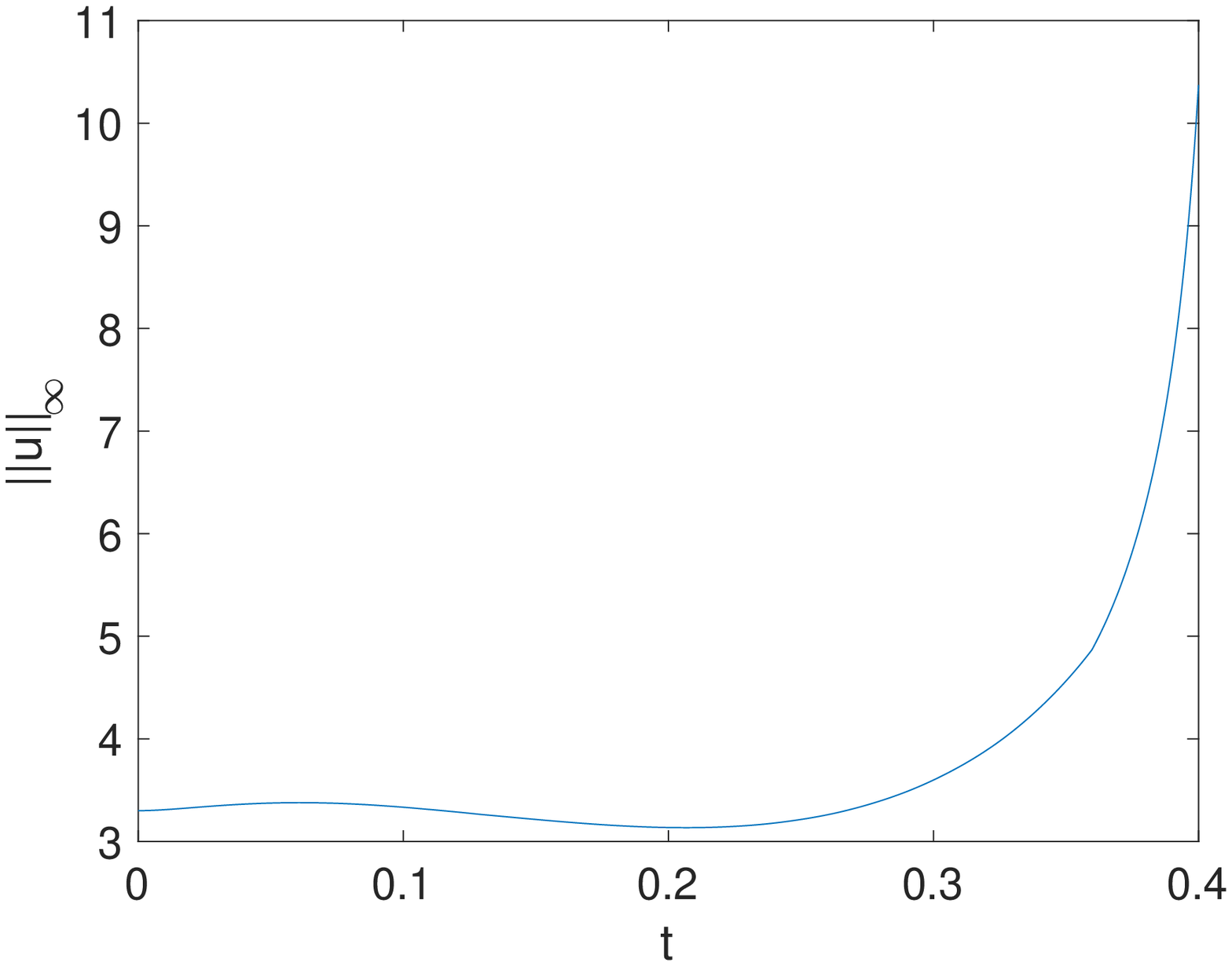}
  \includegraphics[width=0.49\textwidth]{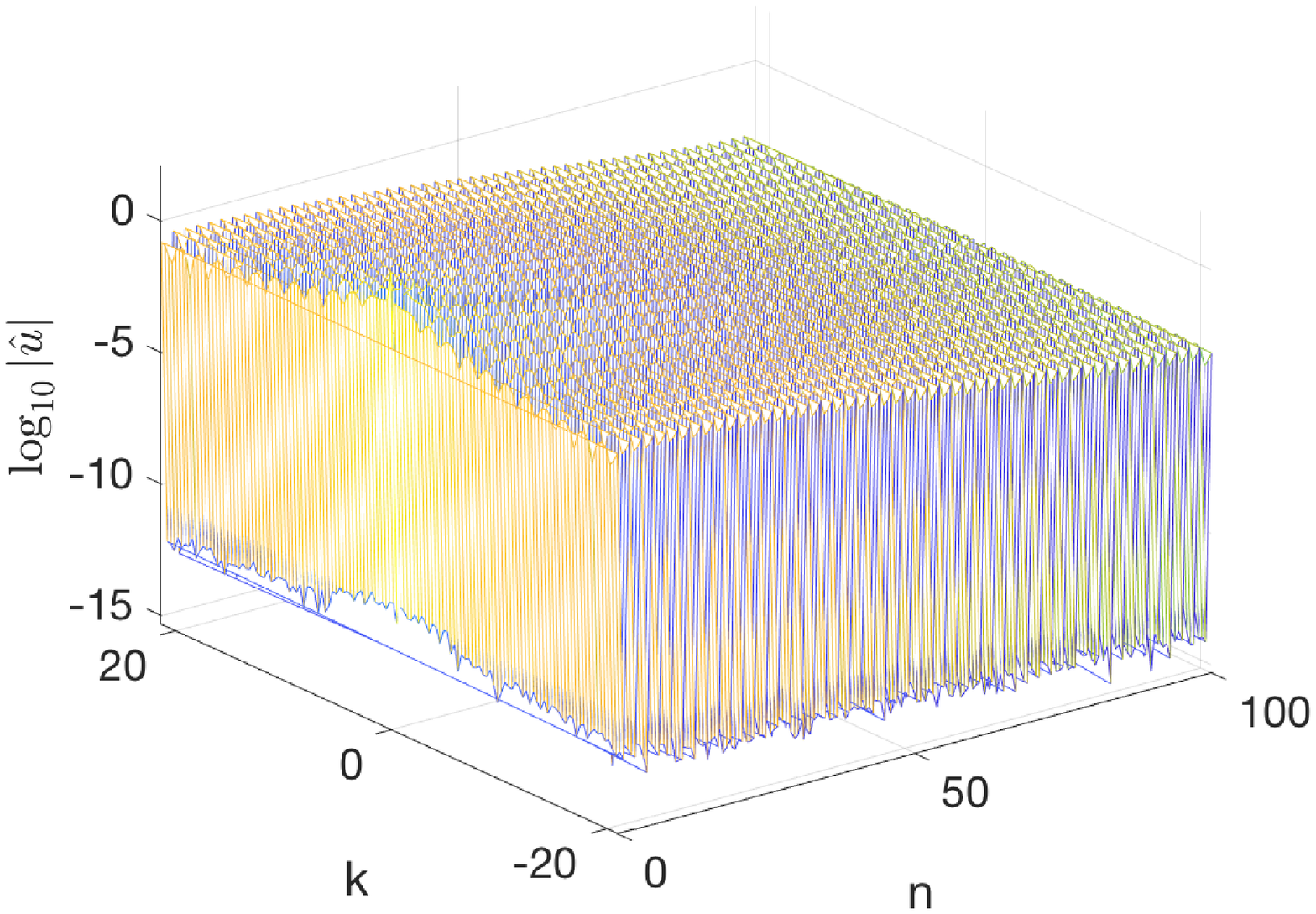}
 \caption{$L^{\infty}$ norm of the solution to the 2D NLS equation for the initial data 
 $(1.1-\exp(-y^{2})) u_{Per}(x,t_{0})$ in dependence of $t$ on the 
 left, and the spectral coefficients of the solution in domain I for 
 the solution at $t=0.4$ on the right.}
 \label{Peregine2d11max}
\end{figure}

A blow-up is suggested also by the $L^{\infty}$ norm of the solution which 
seemingly diverges as can be seen on the left of 
Fig.~\ref{Peregine2d11max}. The spectral coefficients on the right of 
the same figure also indicate a loss of resolution because of said blow-up.

\begin{figure}[htb!]
  \includegraphics[width=0.49\textwidth]{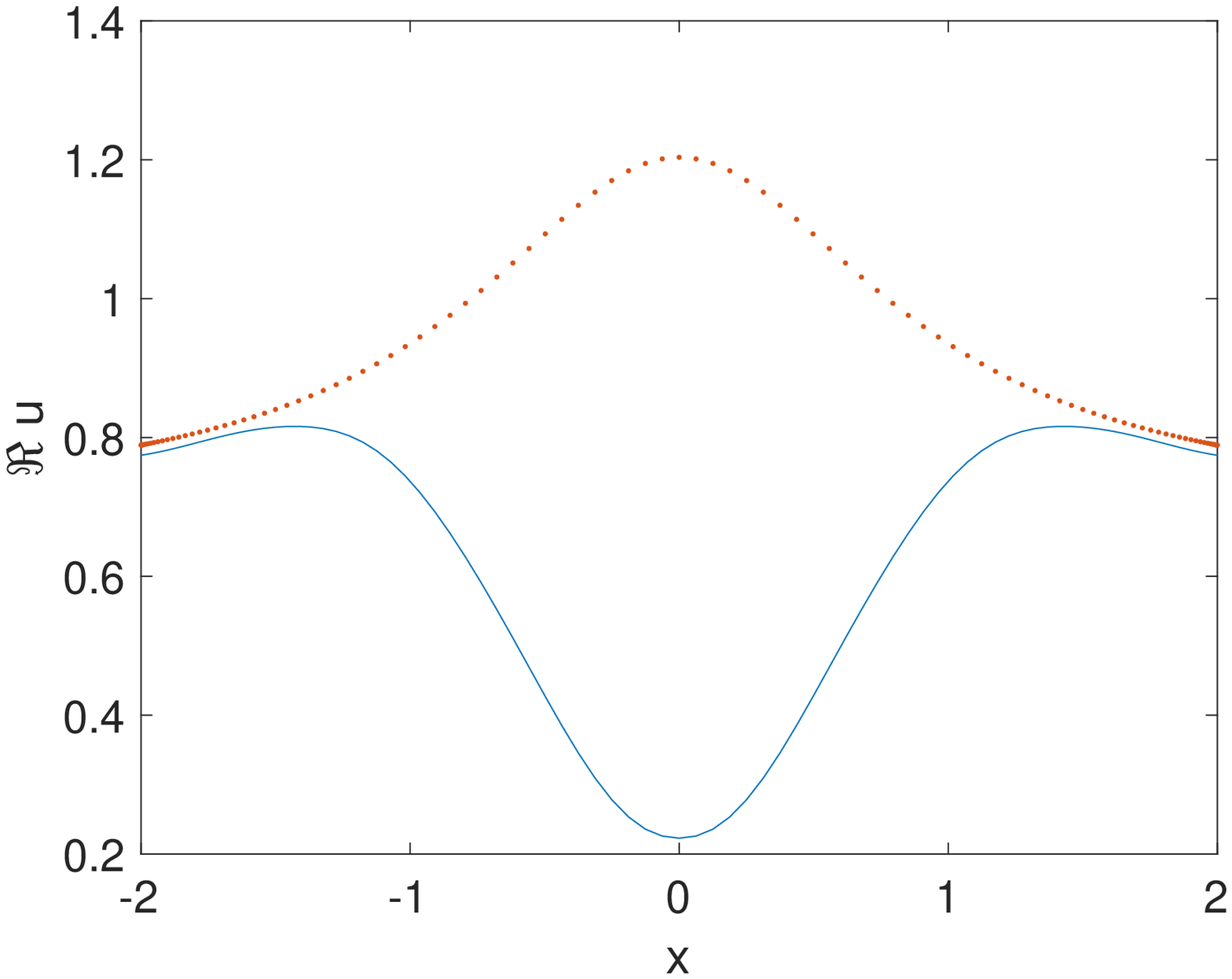}
  \includegraphics[width=0.49\textwidth]{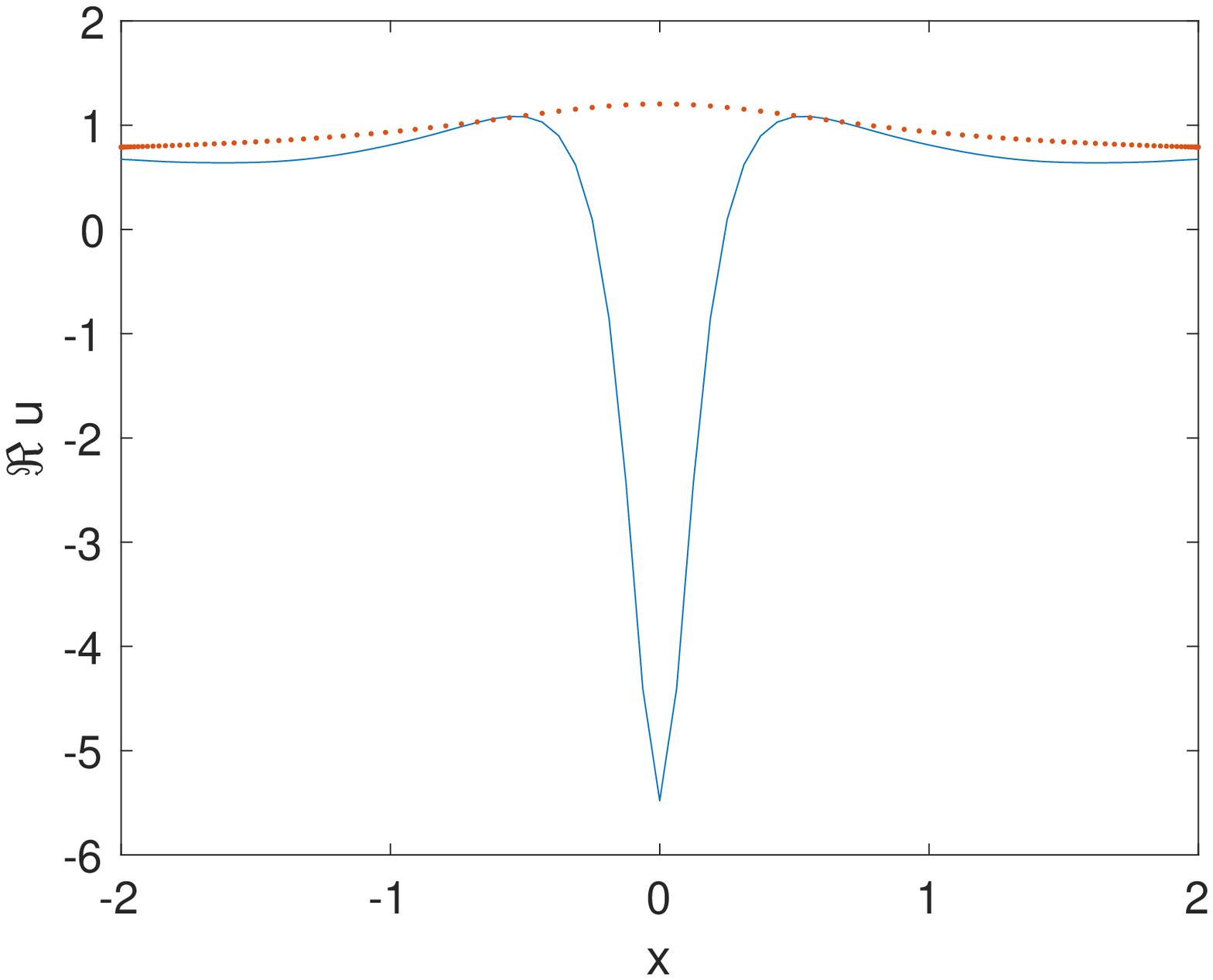}
 \caption{Real part of the solution to the 2D NLS equation for the initial data 
 $(1.1-\exp(-y^{2})) u_{Per}(x,t_{0})$ in blue for $y=0$ the 
 left, and for $y=1.7671$ on the right, and the respective Peregrine 
 solution in red.}
 \label{Peregine2d11y}
\end{figure}

A blow-up once more suggests that the Peregrine solution is strongly 
unstable as a solution to the 2D focusing elliptic NLS equation. Whereas this 
is numerically difficult to decide, the solution is clearly unstable 
in the sense that it does not stay close to the exact solution at the 
same time as can be seen in Fig.~\ref{Peregine2d11y}, where the 
solution for two values of $y$ corresponding to the shown minimum 
$y=0$ and the maximum $y\sim 1.7671$ of $|u|$ in blue and the corresponding 
Peregrine solution in red is given. 

\section{Hyperbolic nonlinear Schr\"odinger equation}
In this section we study some of the examples of the previous 
sections in the context of the hyperbolic NLS equation. Since the 
Peregrine solution is again unstable against all considered 
perturbations, we concentrate on the localized perturbations 
(\ref{initial}) here, and 
show only the solution at the final time. No indication of blow-up is 
found in this case. 

As detailed 
for instance in \cite{Lannes}, hyperbolic NLS or Davey-Stewartson 
equations appear in the modulational regime of the water waves. Since 
the hyperbolic NLS equation has one focusing and one defocusing 
direction, no blow-up is to be expected in this case. Numerically 
this was confirmed in \cite{KS15}, for analytic results see 
\cite{totz}. The defocusing character of the $y$-direction implies 
that the numerical resolution is higher for the examples here (the 
Fourier coefficients in $y$ decrease to lower values), though 
we use the same numerical parameters as in section 3, where also 
plots of the initial configurations can be found.

\begin{figure}[htb!]
 \includegraphics[width=0.49\textwidth]{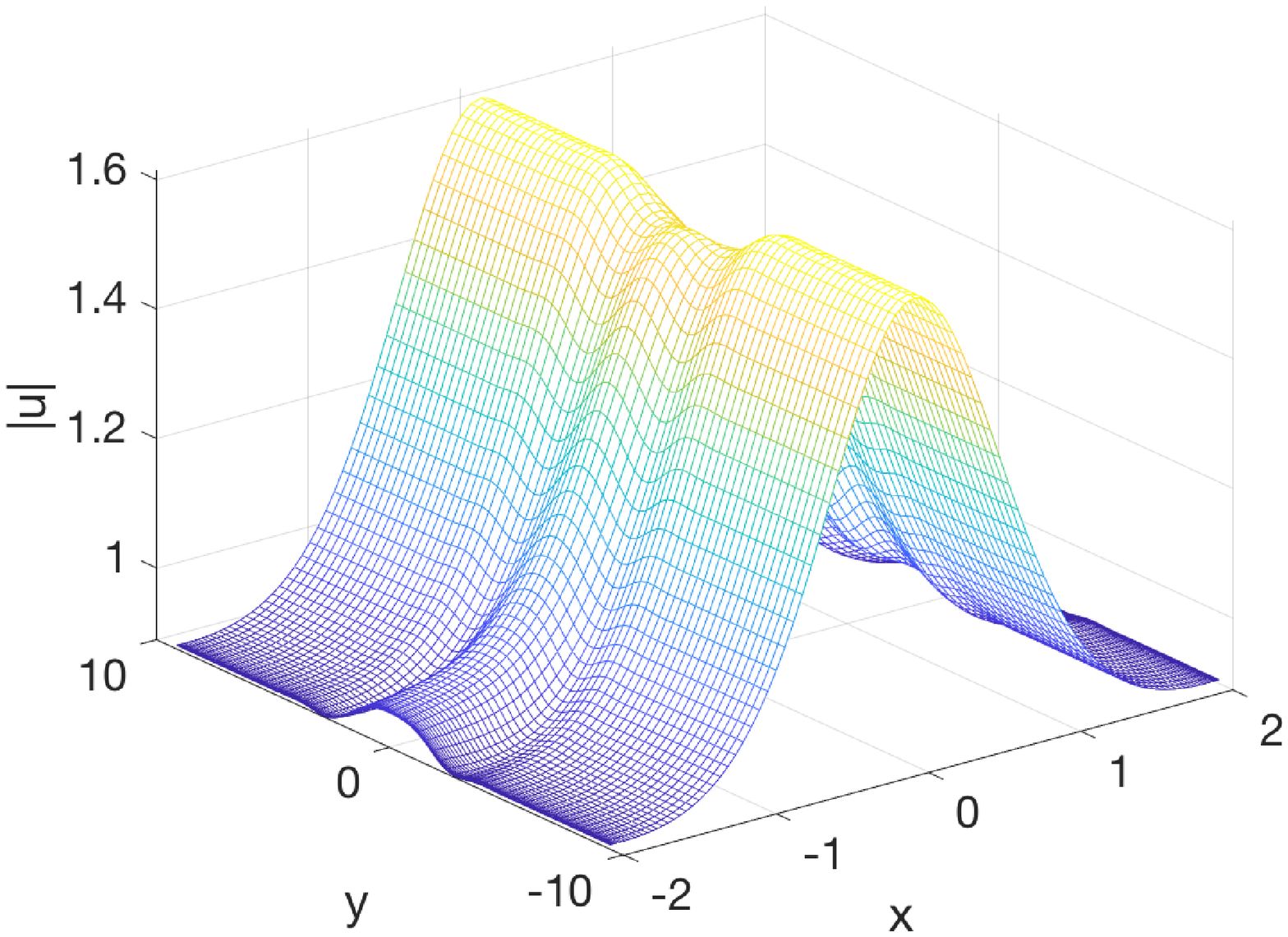}
   \includegraphics[width=0.49\textwidth]{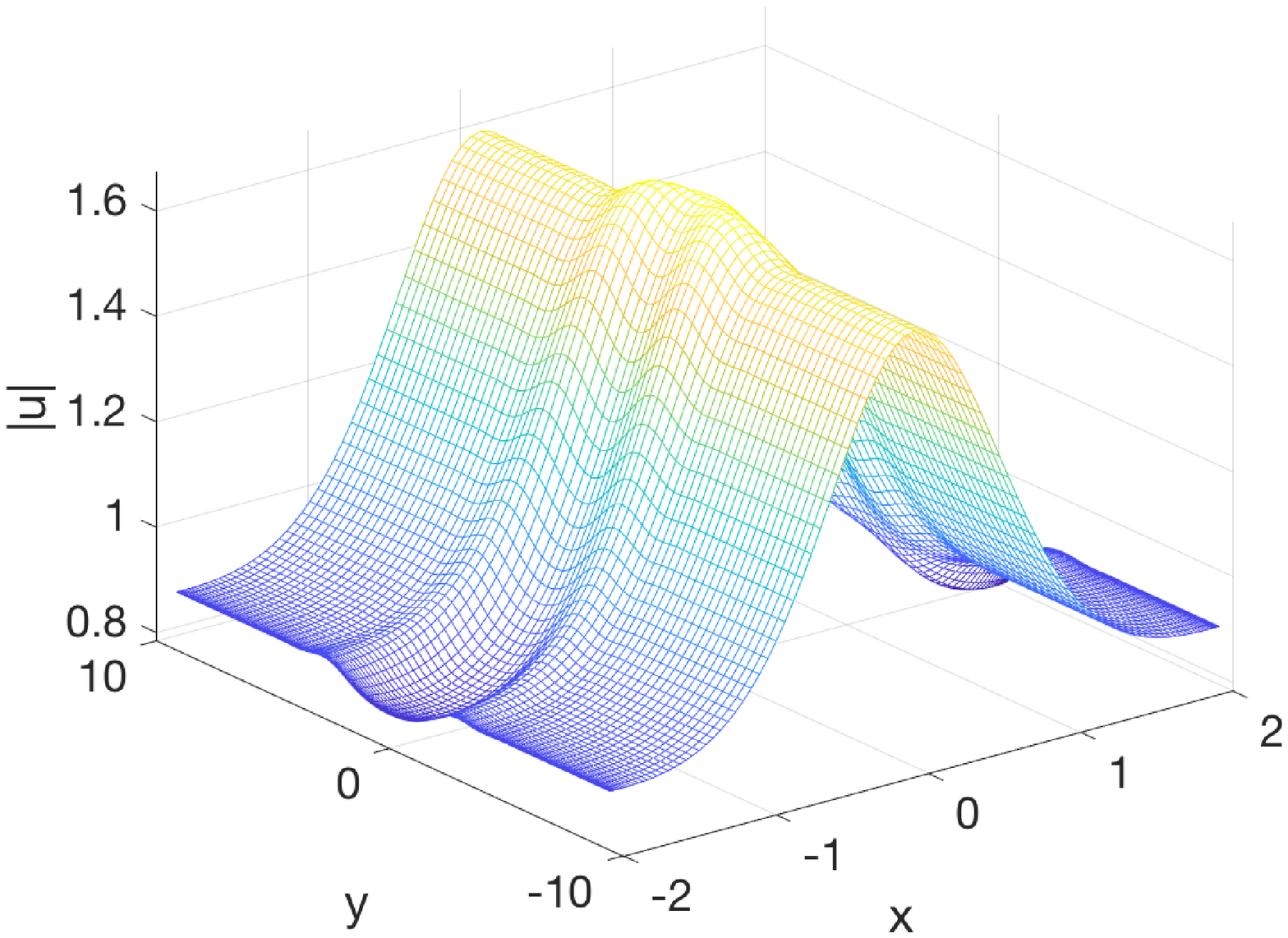}
 \caption{Solution to the hyperbolic NLS equation for the initial data 
 $u(x,y,0)=u_{Per}(x,t_{0})\pm 0.1\exp(-(x+1)^{2}-y^{2})$ for $t=0.5$, on the 
 left for the $+$ sign, on the right for the $-$ sign.}
 \label{peregrine2dhypgaussm1}
\end{figure}

As in Fig.~\ref{Peregine2dp01gaussx1} for the elliptic NLS equation, 
we first consider the case $x_{0}=-1$ and $c=\pm0.1$ in (\ref{initial}). 
The solutions for $t=0.5$ are shown in
Fig.~\ref{peregrine2dhypgaussm1}. Again the perturbations do not disappear.

\begin{figure}[htb!]
 \includegraphics[width=0.49\textwidth]{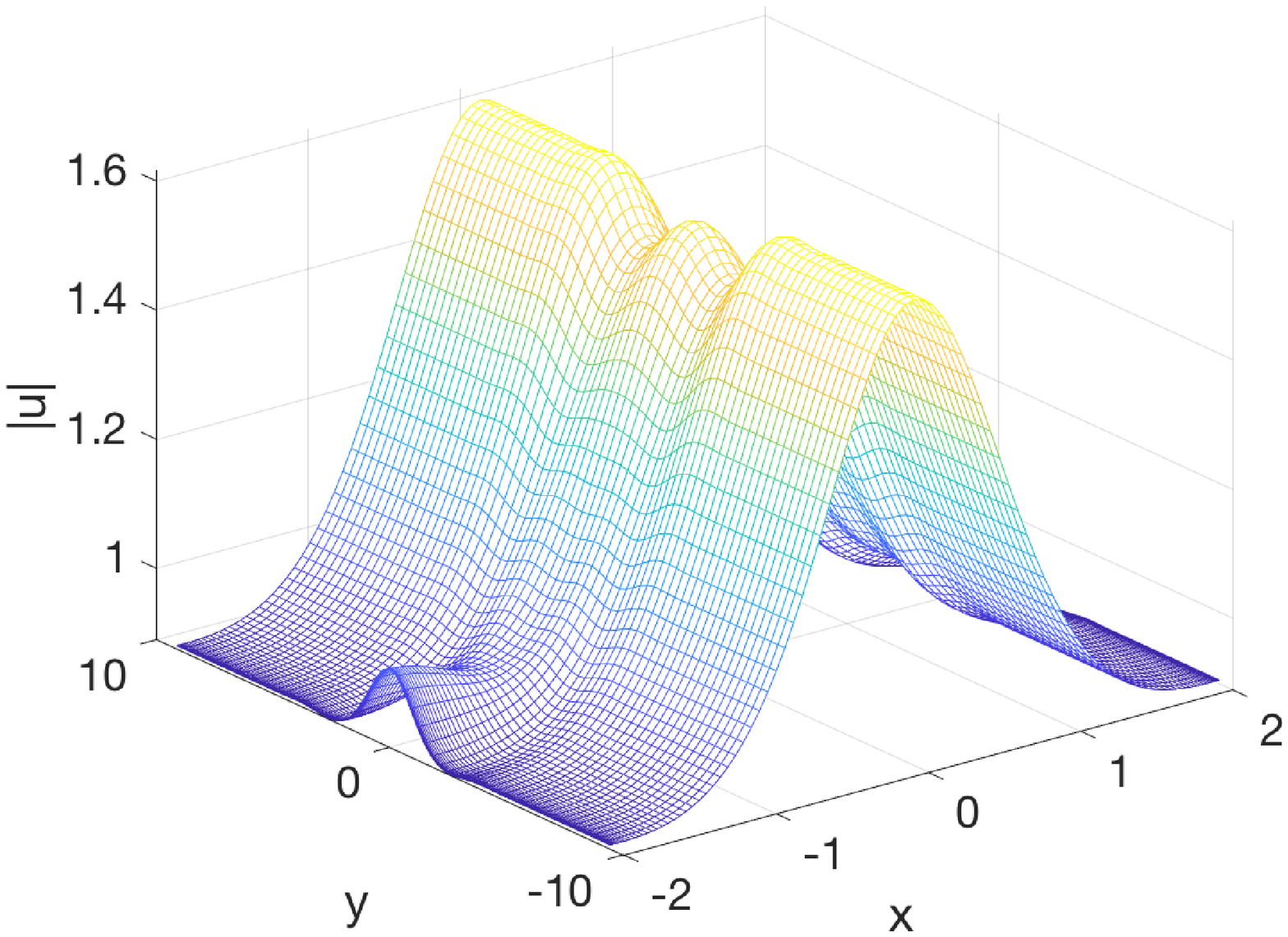}
   \includegraphics[width=0.49\textwidth]{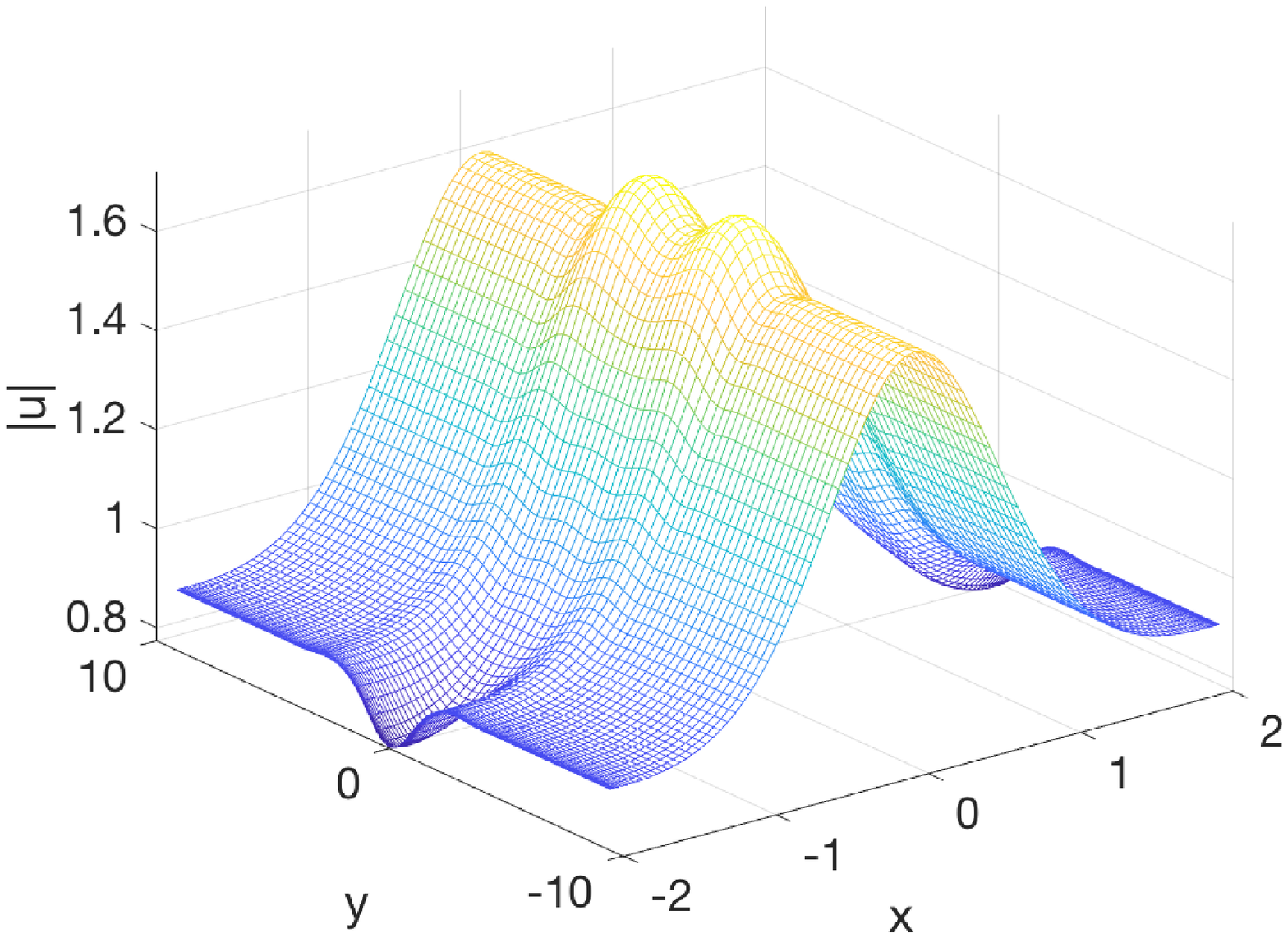}
 \caption{Solution to the hyperbolic NLS equation for the initial data 
 $u(x,y,0)=u_{Per}(x,t_{0})\pm 0.1\exp(-x^{2}-y^{2})$ for $t=0.5$, on the 
 left for the $+$ sign, on the right for the $-$ sign.}
 \label{peregrine2dhypgaussm}
\end{figure}

Of special interest is the situation of the initial data (\ref{initial}) for $x_{0}=0$ and 
$c=\pm0.1$ since a blow-up was conjectured in the case with the $-$ 
sign, see Fig.~\ref{Peregine2dm01gauss}. As can be seen in Fig.~\ref{peregrine2dhypgaussm}, there is no 
indication of a blow-up in this case for either sign of the initial 
data. The solution appears to be simply dispersed, but the initial 
perturbations do not decrease.

\section{Conclusion}
In this paper we have presented a numerical approach to the 2D 
cubic NLS equation (both in the elliptic and hyperbolic case) based on a multi-domain spectral approach 
with a compactified exterior domain in the $x$-coordinate (which 
allows to treat infinity as a finite  point) and a Fourier 
spectral method in the coordinate $y$. This provides a spectral 
approach in both spatial coordinates which allows to achieve a 
numerical error exponentially decreasing with the resolution for 
functions periodic or rapidly decreasing in $y$ and rational and 
bounded at infinity in $x$. For the time integration a fully explicit 
fourth order method was presented based on a 4th order splitting 
scheme and an IRK4 method for the linear step.

This code allows to treat perturbations of the Peregrine solution 
with high precision and efficiency. In the  examples studied the known 
instability of the Peregrine solution in 1D is shown to be amplified 
in 2D, also in the hyperbolic case which has one defocusing direction. All considered perturbations grow in time. In addition,  
numerical evidence for strong instability of the 
Peregrine solution in the elliptic case is presented: the solution appears 
to blow up in finite time if the initial data have an $E[u]$ 
(\ref{energy}) (which is 
finite on the considered torus in $y$) smaller than the one of the 
Peregrine solution.

It is left for a future project to develop better time integration 
schemes for the proposed spatial discretisation. An interesting 
approach would be to use exponential integrators as in 
\cite{KT,HO,etna} for which 
an efficient computation of matrix exponentials will be needed. In 
particular it will be necessary to filter the unphysical eigenvalues 
of the Chebyshev differentiation matrices, see the discussion in 
\cite{trefethen}.


\begin{thebibliography}{99}
    \bibitem{Abl}M.J. Ablowitz and T.P. Horikis, Interacting 
    nonlinear wave envelopes and rogue wave formation in deep water, 
    Physics of Fluids 27, 012107 (2015); https://doi.org/10.1063/1.4906770
 
 \bibitem{bailung} H. Bailung, S. K. Sharma, and Y. Nakamura, Observation of Peregrine solitons in a multicomponent plasma with negative ions, Phys. Rev. Lett. 107, 255005 (2011).
 
 \bibitem{Berenger}
 J. B\'erenger. A perfectly matched layer for the absorption of 
 electromagnetic waves. J. Comput. Phys. 114, 185-200 (1994).



\bibitem{BK} M.~Birem and C.~Klein, Multidomain spectral method 
     for Schr\"odinger equations,  Adv. 
     Comp. Math. DOI 10.1007/s10444-015-9429-9 (2015)

\bibitem{DBU}J.~Bona and J.-C.~Saut, Dispersive Blow-Up II. Schr\"odinger-Type 
Equations, Optical and Oceanic Rogue Waves, Chinese Annals of Math. Series B,
31, (6), 793-810 (2010).

\bibitem{BPSS}J.L. Bona, G. Ponce, J.-C. Saut, C. Sparber, Dispersive 
blow-up for nonlinear Schršdinger equations revisited, JMPA, 102, 
782--811 (2014).

\bibitem{BFS}Boutet de Monvel, A., Fokas, A.S., Shepelsky, D.: Analysis of 
the global relation for the nonlinear Schro ?dinger equation on the 
half-line. Lett. Math. Phys. 65, 199Ð212 (2003)    


\bibitem{cha1}A. Chabchoub, N. P. Hoffmann, and N. Akhmediev, Rogue wave observation in a water wave tank, Phys. Rev. Lett., 106, 204502 (2011).

\bibitem{cha2} A. Chabchoub, N. Hoffmann, M. Onorato, and N. Akhmediev, Super rogue waves: observation of a higher-order breather in water waves, Phys. Rev. X 2, 011015 (2012).

\bibitem{cha3}A. Chabchoub, N. Hoffmann, H. Branger, C. Kharif, and 
N. Akhmediev, Experiments on wind-perturbed rogue wave hydrodynamics 
using the Peregrine breather model, Physics of Fluids, 25, DOI: 10.1063/1.4824706 (2013).

\bibitem{cha4}
A. Chabchoub, K. Mozumi, N. Hoffmann, A.V. Babanin, 
A. Toffoli, J.N. Steer, T.S. van den Bremer, N. 
Akhmediev, M. Onorato, and T. Waseda, Directional soliton and breather beams   
PNAS, 116 (20) 9759-9763, 2019.


\bibitem{dubard} P.~Dubard,  P.~Gaillard,  C.~Klein, and V.B.~Matveev,   On multi-rogue wave solutions of the NLS 
equation and positon solutions of the KdV equation, Eur. Phys. J. 
Special Topics,  185, 247--258 (2010).	


%



\bibitem{GO} Grosch, C.E., Orszag, S.A.: Numerical solution of problems in unbounded regions: coordinate
transforms. J. Comput. Phys. 25, 273Ð296 (1977)

\bibitem{HO}M. Hochbruck and A. Ostermann, \emph{Exponential 
integrators}, Acta Numerica (2010), pp. 209Ð286.

\bibitem{KT}A.-K. Kassam and L. Trefethen, Fourth-Order Time-Stepping for stiff PDEs, SIAM J. Sci.
Comput., 26 (2005), pp. 1214Ð1233.

\bibitem{Ka}T. Kato, Trotter's product formula for an arbitrary pair of self-adjoint contraction semigroups,
vol. 3, Academic Press, Boston, 1978, pp. 185Ð195.

\bibitem{KBAA}U. Al Khawaja, H. Bahlouli, M. Asad-uz-zaman, S.M. 
Al-Marzoug, Modulational instability analysis of the Peregrine 
soliton, Commun Nonlinear Sci Numer Simulat 19, 2706Ð2714 (2014).

\bibitem{kibler}B. Kibler, J. Fatome, C. Finot, G. Millot, F. Dias, 
G. Genty, N. Akhmediev, and J. M. Dudley, The Peregrine soliton in 
nonlinear fibre optics, Nat. Phys. 6, 790--795 (2010). 


\bibitem{etna} C.~Klein,  
Fourth order time-stepping for low dispersion Korteweg-de Vries and nonlinear Schr\"odinger equation,  
ETNA, 29, 116--135 (2008).

\bibitem{KH} C. Klein and M. Haragus, Numerical study of the stability of the Peregrine breather, Annals of
Mathematical Sciences and Applications, 2(2), 217-239 (2017).

\bibitem{KS15} C.~Klein and J.-C.~Saut, \emph{A numerical approach to 
     Blow-up issues for Davey-Stewartson II type systems},  Comm. 
	 Pure Appl. Anal.  14:4,  1443--1467     (2015) 
     
     


\bibitem{KS}C. Klein and N. Stoilov, A numerical study of blow-up 
mechanisms for Davey-Stewartson II systems, Stud. Appl. Math., DOI : 
10.1111/sapm.12214 (2018)  

\bibitem{Lannes}D. Lannes, Water waves: mathematical theory and 
asymptotics, Mathematical Surveys and Monographs, vol 188 (2013), AMS, Providence.     

\bibitem{Li}Y.C. Li, On the so-called rogue waves in the nonlinear 
Schr\"odinger equation, \emph{preprint}, arXiv:1511.00620v1 (2015)

\bibitem{MR}F. Merle and P. Raphael. On universality of blow-up 
profile for l2 critical nonlinear Schr\"odinger equation. Inventiones
Mathematicae, 156:565Ð672, 2004.

\bibitem{munoz}	C. Mu\~noz, Instability in nonlinear Schr\"odinger breathers, 
Proyecciones Vol. 36, no. 4 (2017) 653Ð683.



\bibitem{Peregrine}D.H.~Peregrine, Water waves, nonlinear 
Schr\"odinger equations and their solutions, J. Austral. Math. Soc. 
B, 25, 16--43, DOI:10.1017/S0334270000003891 (1983).


\bibitem{sulem}C. Sulem and P.L. Sulem. The nonlinear Schr\"odinger equation. Springer, 1999.

\bibitem{totz}N. Totz, Global well-posedness of 2D non-focusing 
Schr\"odinger equations via
rigorous modulation approximation, J. Dif. Equat. 261:2251Ð2299 
(2016). 


\bibitem{trefethen} 
L.~N.~Trefethen,  
\emph{ Spectral Methods in Matlab}, SIAM, Philadelphia, PA, 2000.  


\bibitem{TK}H. Trotter, On the Product of Semi-Groups of Operators, Proceedings of the American
Mathematical Society, 10 (1959), pp. 545Ð551.
\bibitem{WR}Weideman, J.A.C. and Reddy, S.C., A Matlab differentiation matrix suite, ACM TOMS, 26
(2000), 465--519.

\bibitem{lorene} www.lorene.obspm.fr
\bibitem{Y}H. Yoshida, Construction of higher Order symplectic Integrators, Physics Letters A, 150
(1990), pp. 262--268.

\bibitem{ZS} V.E. Zakharov and A.B.  Shabat, Exact theory of two-dimensional self-focusing and one-dimensional selfmodulation of waves in nonlinear media. Sov. Phys. JETP, 34(1), 62--69 
(1972); translated from Zh. Eksp. Teor. Fiz. 1, 118--134 (1971).

\bibitem{Zhengpml}
C. Zheng. A perfectly matched layer approach to the nonlinear 
Schr\"odinger wave equations. J. Comput. Phys. 227, 537-556 (2007).

\bibitem{Zhengtbc}
C. Zheng, Exact nonreflecting boundary conditions for one-dimensional 
cubic nonlinear Schr\"odinger equations, J. Comput. Phys. 
215, 552-565 (2006).

\end{thebibliography}
\end{document}